\pdfoutput=1
\RequirePackage{ifpdf}
\ifpdf 
\documentclass[pdftex]{sigma}
\else
\documentclass{sigma}
\fi

\newcommand{\Dchainthree}[6]{
\rule[-3\unitlength]{0pt}{8\unitlength}
\begin{picture}(26,5)(0,3)
\put(1,2){\ifthenelse{\equal{#1}{l}}{\circle*{2}}{\circle{2}}}
\put(2,2){\line(1,0){10}}
\put(13,2){\ifthenelse{\equal{#1}{m}}{\circle*{2}}{\circle{2}}}
\put(14,2){\line(1,0){10}}
\put(25,2){\ifthenelse{\equal{#1}{r}}{\circle*{2}}{\circle{2}}}
\put(1,5){\makebox[0pt]{\scriptsize #2}}
\put(7,4){\makebox[0pt]{\scriptsize #3}}
\put(13,5){\makebox[0pt]{\scriptsize #4}}
\put(19,4){\makebox[0pt]{\scriptsize #5}}
\put(25,5){\makebox[0pt]{\scriptsize #6}}
\end{picture}}

\newcommand{\Dchaintwo}[4]{
\rule[-3\unitlength]{0pt}{8\unitlength}
\begin{picture}(14,5)(0,3)
\put(1,2){\ifthenelse{\equal{#1}{l}}{\circle*{2}}{\circle{2}}}
\put(2,2){\line(1,0){10}}
\put(13,2){\ifthenelse{\equal{#1}{r}}{\circle*{2}}{\circle{2}}}
\put(1,5){\makebox[0pt]{\scriptsize #2}}
\put(7,4){\makebox[0pt]{\scriptsize #3}}
\put(13,5){\makebox[0pt]{\scriptsize #4}}
\end{picture}}

\newcommand{\lact }{.}

\newcommand{\btxandshort}[1]{and} \newcommand{\btxpagesshort}[1]{pp.} \newcommand{\Btxinshort}[1]{In}
\newcommand{\btxphdthesis}[1]{phd-thesis} \newcommand{\btxeditorshort}[1]{Ed.} \newcommand{\btxeditorsshort}[1]{Eds.}
\newcommand{\btxvolumeshort}[1]{vol.} \newcommand{\btxofseriesshort}[1]{ser.}

\newcommand{\SLZ}[1]{\mathrm{SL}(#1,\mathbb{Z})}

\newcommand{\rersys}[1]{\boldsymbol{\Delta}{}^{#1 \mathrm{re}}}

\numberwithin{table}{section}
\numberwithin{equation}{section}
\newtheorem{Theorem}{Theorem}[section]
\newtheorem{Corollary}[Theorem]{Corollary}
\newtheorem{Lemma}[Theorem]{Lemma}
\newtheorem{Proposition}[Theorem]{Proposition}
{ \theoremstyle{definition}
\newtheorem{Definition}[Theorem]{Definition}
\newtheorem{Remark}[Theorem]{Remark} }

\newlength{\mpb}

\DeclareMathOperator{\ad}{ad}
\DeclareMathOperator{\Aut}{Aut}
\DeclareMathOperator{\Hom}{Hom}
\DeclareMathOperator{\End}{End}

\begin{document}

\newcommand{\arXivNumber}{1407.6817}

\allowdisplaybreaks

\renewcommand{\PaperNumber}{011}

\FirstPageHeading

\ShortArticleName{Rank 2 Nichols Algebras of Diagonal Type over Fields of Positive Characteristic}

\ArticleName{Rank 2 Nichols Algebras of Diagonal Type\\
over Fields of Positive Characteristic}

\Author{Jing WANG and Istv{\'a}n HECKENBERGER}
\AuthorNameForHeading{J.~Wang and I.~Heckenberger}
\Address{Philipps-Universit\"at Marburg, FB Mathematik und Informatik,\\
Hans-Meerwein-Stra{\ss}e, 35032 Marburg, Germany}
\Email{\href{mailto:jing@mathematik.uni-marburg.de}{jing@mathematik.uni-marburg.de},
\href{mailto:heckenberger@mathematik.uni-marburg.de}{heckenberger@mathematik.uni-marburg.de}}
\URLaddress{\url{http://www.mathematik.uni-marburg.de/~jing/},\\
\hspace*{10.5mm}\url{http://www.mathematik.uni-marburg.de/~heckenberger/}}

\ArticleDates{Received August 01, 2014, in f\/inal form February 02, 2015; Published online February 07, 2015}

\Abstract{The paper introduces a~new method to determine all rank two Nichols algebras of diagonal type over f\/ields of
positive characteristic.}


\Classification{16T05}

\section{Introduction}
The theory of Nichols algebras was dominated and motivated by Hopf algebra theory.
In 1978, W.~Nichols f\/irst introduced the structure of Nichols algebra in the paper ``Bialgebras of type
one''~\cite{n-78}, where he studied certain pointed Hopf algebras.
S.L.~Woronowicz rediscovered this structure in his approach to ``quantum dif\/ferential calculus''~\cite{w1987,w1989}.
M.~Rosso and G.~Lusztig def\/ined and used them to present quantum groups in a~dif\/ferent language~\cite{l-2010, Rosso98}.
In fact, Nichols algebra turns out to be very important in Hopf algebras and quantum
groups~\cite{a-AndrGr99,inp-AndrSchn02,a-Schauen96} and has applications in conformal f\/ield theory and mathematical
physics~\cite{Semi-2011,Semi-2012,Semi-2013}.
Nichols algebras play an important role in the classif\/ication of pointed Hopf algebras with certain f\/initeness
properties in papers by N.~Andruskiewitsch and H.-J.~Schneider, see for example~\cite{a-AndrSchn98,inp-AndrSchn02}.
The crucial step to classify pointed Hopf algebras is the computation of the Nichols algebras.
The explicit presentations by generators and relations of a~f\/inite-dimensional Nichols algebra of a~braided vector space
in suitable classes are crucial for the classif\/ication lifting method in~\cite{a-AndrSchn98,inp-AndrSchn02}.

Several authors have already classif\/ied both inf\/inite and f\/inite dimensional Nichols algebra of Cartan type,
see~\cite{a-AndrSchn00,a-Heck06a,Rosso98}.
Further, N.~Andruskiewitsch~\cite{inp-Andr02} stated the following question.

\textbf{Question 5.9.} Given a~braiding matrix $(q_{ij})_{1\leq i,j\leq \theta} $ whose entries are roots of~$1$, when
$\mathcal{B}(V)$ is f\/inite-dimensional, where~$V$ is a~vector space with basis $x_1,\dots,x_{\theta}$ and braiding $c(x_i\otimes
x_j)=q_{ij}(x_j\otimes x_i)$? If so, compute $\dim_\Bbbk \mathcal{B}(V)$, and give a~``nice'' presentation by generators and
relations.

The f\/irst half of Question 5.9 was answered by the f\/irst named author in~\cite{a-Heck09} when the characteristic of the
f\/ield is~0.
The crucial theoretical tools of the classif\/ication were the Weyl groupoid of a~braided vector space of diagonal type
and the root system associated to a~Nichols algebra of diagonal type, see~\cite{a-Heck06a}.
From V.~Kharchenko~\cite[Theorem~2]{a-Khar99} any Nichols algebra~$\mathcal{B}(V)$ of diagonal type has a~(restricted) Poincar\'e--Birkhof\/f--Witt
basis consisting of homogenous elements with respect to the $\mathbb{Z}^n$-grading of~$\mathcal{B}(V)$.
In~\cite{a-Heck06a}, the root system and the Weyl groupoid of $\mathcal{B}(V)$ for a~Nichols algebra $\mathcal{B}(V)$ of diagonal type
was def\/ined.
This Weyl groupoid plays a~similar role as the Weyl group does for ordinary root systems.
Based on these results, in~\cite{HY08} and~\cite{c-Heck09b} the abstract combinatorial theory of Weyl groupoids and
generalized root systems was initiated.
Later, the theory of root systems and Weyl groupoids was extended to more general Nichols algebras
in~\cite{a-AHS08,HS10,a-HeckSchn12a}.
With the classif\/ication result in~\cite{a-Heck09}, N.~Andruskiewitsch and \mbox{H.-J.}~Schnei\-der~\cite{a-AndrSchn05} obtained
a~classif\/ication theorem about f\/inite-dimensional pointed Hopf algebras under some technical assumptions.
On the other hand, it is natural and desirable to analyze the classif\/ication of Nichols algebras of diagonal type for
arbitrary f\/ields.
The authors in~\cite{clw} constructed new examples of Nichols algebras in positive characteristic by applying
a~combinatorial formula for products in Hopf quiver algebras.

In this paper, all rank 2 Nichols algebras of diagonal type with a~f\/inite root system over f\/ields of positive
characteristic are classif\/ied.
We introduce some properties of rank two Cartan graphs, see Theorem~\ref{theo.cartan}.
Theorem~\ref{theo.cartan} characterizes f\/inite connected Cartan graphs of rank two in terms of certain integer
sequences.
This theorem simplif\/ies substantially the calculations needed to check that the Weyl groupoids of the Nichols algebras
in Tables~\ref{tab.1}--\ref{tab.6} of our classif\/ication are f\/inite.
Indeed, using the theorem, these calculations can be done by hand within a~very short time, in contrast to the
calculations based on the def\/inition of the Weyl groupoid $\mathcal{W}$ or using a~longest element of $\mathcal{W}$.
The main result of this paper is Theorem~\ref{Theo:clasi}.
Table~7 illustrates all the exchange graphs of the corresponding Cartan graphs in Theorem~\ref{Theo:clasi}.

The structure of the paper is as follows.
In Section~\ref{se:Cschemes} we recall the def\/inition of Cartan graphs, their Weyl groupoids and root systems.
Some well-known results are also recalled.
In Section~\ref{se:Car of Nic}, Theorem~\ref{theo.regualrcar} associates a~semi-Cartan graph $\mathcal{C}(M)$ of rank~$\theta$
to a~tuple~$M$ of f\/inite-dimensional irreducible Yetter--Drinfel'd-modules.
Further, $\mathcal{C}(M)$ is a~Cartan graph if the set of real roots of~$M$ is f\/inite.
The Dynkin diagram for a~braided vector space of diagonal type is recalled in Section~\ref{se:dynkin diag} and some
corollaries are also obtained, see Lemma~\ref{le:Dynkin} and Proposition~\ref{prop:Cs}.
In Section~\ref{se:classi of Carscheme} we recall certain integer sequences and prove a~local property of them in
Theorem~\ref{Theo:clascar}.
This is another main technical part for the classif\/ication.
Finally the main result of this paper is formulated in Section~\ref{se:cla_Nichol}, see Theorem~\ref{Theo:clasi}.
In this section, all rank two Nichols algebras of diagonal type with a~f\/inite root system over f\/ields of positive
characteristic are classif\/ied.
Since many subcases have to be considered, this is the largest part of the paper.
To simplify the results, this paper ends with 5 tables containing all the Dynkin diagrams of rank two braided vector
spaces with a~f\/inite root system for all f\/ields of positive characteristic.

Throughout the paper $\Bbbk$ denotes a~f\/ield of characteristic $p> 0$.
Let $\Bbbk^*=\Bbbk\setminus \{0\}$.
The set of natural numbers not including $0$ is denoted by $\mathbb{N}$ and we write $\mathbb{N}_0=\mathbb{N}\cup \{0\}$.
For $n\in \mathbb{N}$, let $G_n'$ denote the set of primitive $n$-th roots of unity in $\Bbbk$, that is $G'_n=\{q\in \Bbbk^*\,|\,
q^n=1$, $q^k\not=1~\text{for all}~1\leq k < n\}$.

\section{Cartan graphs and root systems}\label{se:Cschemes}

Let~$I$ be a~non-empty f\/inite set.
Recall from~\cite[\S~1.1]{b-Kac90} that a~generalized Cartan matrix is a~matrix $A=(a_{ij})_{i,j\in I}$ with integer entries such that
\begin{itemize}\itemsep=0pt
\item $a_{ii}=2$ and $a_{jk}\le 0$ for any $i,j,k\in I$ with $j\not=k$,
\item if $a_{ij}=0$ for some $i,j\in I$, then $a_{ji}=0$.
\end{itemize}
Let $\mathcal{X}$ be a~non-empty set and let $r: I\times \mathcal{X}\to \mathcal{X}$ be a~map.
For all $i\in I$, let $r_i: \mathcal{X}\to \mathcal{X}$, $ X \mapsto r(i,X)$.
Let $A^X =(a_{ij}^X)_{i,j \in I}$ be a~generalized Cartan matrix in $\mathbb{Z}^{I \times I}$ for all $X\in \mathcal{X}$.
The quadruple $\mathcal{C}= \mathcal{C}(I,\mathcal{X},r,(A^X)_{X \in \mathcal{X}})$ is called a~{\textit{semi-Cartan graph}} if
\begin{itemize}\itemsep=0pt
\item $r_i^2 = \mathrm{id}_{\mathcal{X}}$ for all $i \in I$,
\item $a^X_{ij} = a^{r_i(X)}_{ij}$ for all $X\in \mathcal{X}$ and $i,j\in I$.
\end{itemize}

The cardinality of~$I$ is called the {\textit{rank}} of $\mathcal{C}$, and the elements of~$I$ the {\textit{labels}} of $\mathcal{C}$.
The elements of $\mathcal{X}$ are called the {\textit{points}} of $\mathcal{C}$.
Semi-Cartan graphs are called Cartan schemes in~\cite{c-Heck09a}.
We change the terminology in order to increase the recognizability of our structures, to avoid possible confusion with
other mathematical concepts, and in order to shorten the notation for our main objects, the Cartan graphs.
We thank N.~Andruskiewitsch and H.-J.~Schneider for a~very fruitful discussion on this issue.

The {\textit{exchange graph}} of $\mathcal{C}$ is a~labeled non-oriented graph with vertices corresponding to points of $\mathcal{C}$,
and edges marked by labels of $\mathcal{C}$, where two vertices $X$, $Y$ are connected by an edge~$i$ if and only if $X\not=Y$ and
$r_i(X)=Y$ (and $r_i(Y)=X$).
For simplif\/ication, instead of several edges we display only one edge with several labels.
A~semi-Cartan graph is called {\textit{connected}} if its exchange graph is a~connected graph.

Let $\mathcal{C}= \mathcal{C}(I,\mathcal{X},r,(A^X)_{X \in \mathcal{X}})$ be a~semi-Cartan graph.
We f\/ix once and for all the notation $(\alpha_i)_{i\in I}$ for the standard basis of $\mathbb{Z}^I$.
Then there exists a~unique category $\mathcal{D}(\mathcal{X}, I)$ with objects Ob$\mathcal{D}(\mathcal{X}, I)=\mathcal{X}$ and morphisms
$\Hom(X,Y)=\{(Y,f,X)\,|\,f\in \End(\mathbb{Z}^I)\}$ for $X, Y\in \mathcal{X}$, such that the composition is def\/ined~by
\begin{gather*}
(Z,g,Y)\circ (Y,f,X)=(Z, gf, X)
\end{gather*}
for all $ X, Y, Z\in \mathcal{X}$, $f,g\in \textrm{End}(\mathbb{Z}^I)$. For all $X\in \mathcal{X}$ and all $i\in I$, let
\begin{gather}
\label{eq-si}
s_i^X\in \Aut\big(\mathbb{Z}^I\big),
\qquad
s_i^X \alpha_j=\alpha_j-a_{ij}^X \alpha_i
\end{gather}
for all $j\in I$.
We write {\textit{$\mathcal{W}(\mathcal{C})$}} for the smallest subcategory of $\mathcal{D}(\mathcal{X}, I)$ which contains all morphisms $(r_i(X),
s_i^X, X)$, where $i\in I $ and $X\in \mathcal{X}$.
The morphisms $(r_i(X), s_i^X, X)$ are usually abbreviated by $s_i^X$, or by $s_i$, if no confusion is possible.
Since all generators are invertible $\mathcal{W}(\mathcal{C})$ is a~groupoid.

Let $\mathcal{C}= \mathcal{C}(I,\mathcal{X},r, (A^X)_{X \in \mathcal{X}})$ be a~semi-Cartan graph and $(\boldsymbol{\Delta}^{X})_{X\in
\mathcal{X}}$ a~family of sets $\boldsymbol{\Delta}^X \subset \mathbb{Z}^I$.
We say that $\mathcal{R}=\mathcal{R}(\mathcal{C}, (\boldsymbol{\Delta}^{X})_{X\in \mathcal{X}})$ is a~{\textit{root system}} of type $\mathcal{C}$
if and only if
\begin{itemize}
\itemsep=0pt
\item $\boldsymbol{\Delta}^X=\big(\boldsymbol{\Delta}^ X\cap \mathbb{N}_0^I\big)\cup -\big(\boldsymbol{\Delta}^X\cap \mathbb{N}_0^I\big)$,
\item $\boldsymbol{\Delta}^X\cap \mathbb{Z}\alpha_i=\{\alpha_i,-\alpha_i\}$ for all $i\in I$,
\item $s_i^X\big(\boldsymbol{\Delta}^X\big)= \boldsymbol{\Delta}^{r_i(X)}$ for all $i \in I$,
\item $(r_i r_j)^{m_{ij}^X}(X) = X$ for any $i,j \in I$ such that $i\not=j$ where $m_{ij}^X:= \big|\boldsymbol{\Delta}^X\cap
(\mathbb{N}_0\alpha_i + \mathbb{N}_0 \alpha_j)\big|$ is f\/inite.
\end{itemize}
For any category $\mathcal{D}$ and any object~$X$ in $\mathcal{D}$, let $\Hom(\mathcal{D}, X)=\cup_{Y\in \mathcal{D}}\Hom(Y,X)$.
For all $X\in \mathcal{X}$, the set
\begin{gather*}
\rersys{X}=\big\{w\alpha_i \in \mathbb{Z}^I\,\big|\,w\in \Hom(\mathcal{W}(\mathcal{C}),X)\big\}
\end{gather*}
is called the {\textit{set of real roots}}of $\mathcal{C}$ at~$X$.
The elements of $\rersys{X}_{\boldsymbol{+}}=\rersys{X}\cap \mathbb{N}_0^I$ are called \textit{positive roots} and those of
$\rersys{X}\cap -\mathbb{N}_0^I$ {\textit{negative roots}}.
A~semi-Cartan graph is called {\textit{finite}} if $\rersys{X}$ is a~f\/inite set for all $X\in \mathcal{X}$.

Let $t_{ij}^X=\big|\rersys{X}\cap (\mathbb{N}_0 \alpha_i+\mathbb{N}_0 \alpha_j)\big|$. We say that $\mathcal{C}$ is a~{\textit{Cartan graph}} if the
following hold:
\begin{itemize}
\itemsep=0pt
\item For all $X\in \mathcal{X}$ the set $\rersys{X}$ consists of positive and negative roots.
\item Let $X\in \mathcal{X}$ and $i, j\in I$.
If $t_{ij}^X<\infty$ then $(r_ir_j)^{t_{ij}^X}(X)=X$.
\end{itemize}
In that case, $\mathcal{W}(\mathcal{C})$ is called the {\textit{Weyl groupoid}} of $\mathcal{C}$.
\begin{Remark}\quad
\begin{itemize}\itemsep=0pt
\item A~semi-Cartan graph $\mathcal{C}$ is a~Cartan graph if and only if $\mathcal{R}=\mathcal{R}\big(\mathcal{C},(\rersys{X})_{X\in \mathcal{X}}\big)$
is a~root system of type $\mathcal{C}$.
\item For any f\/inite Cartan graph $\mathcal{C}$, there is a~unique root system $\mathcal{R}=\mathcal{R}\big(\mathcal{C},(\rersys{X})_{X\in \mathcal{X}}\big)$
of type $\mathcal{C}$, see~\cite[Propositions~2.9, 2.12]{c-Heck09b}.
\end{itemize}
\end{Remark}

\section{Cartan graphs for Nichols algebras}
\label{se:Car of Nic}
Let~$H$ be a~Hopf algebra over $\Bbbk$ with bijective antipode.
Let $ {}^H_H\mathcal{YD}$ denote the category of Yetter--Drinfel'd modules over~$H$ and $\mathcal{F}_\theta^H $ the set of
$\theta$-tuples of f\/inite-dimensional irreducible objects in ${}^H_H\mathcal{YD}$ for all $\theta \in \mathbb{N}$.
For all $V\in {}^H_H\mathcal{YD}$, let~$\delta$ denote the left coaction of~$H$ on~$V$ and $\cdot$
the left action of~$H$ on~$V$.
Let $\theta \in \mathbb{N}$, $I=\{1,\dots, \theta\}$, and $M = (M_1,\dots,M_{\theta})\in \mathcal{F}_\theta^H$.
Write $[M]=([M_1], \dots, [M_{\theta}])\in \mathcal{X}_\theta^H$, where $\mathcal{X}_\theta^H $ denotes the set of $\theta $-tuples of isomorphism
classes of f\/inite-dimensional irreducible objects in $ {}^H_H\mathcal{YD}$.
Let $\mathcal{B}(M)$ denote the Nichols algebra $\mathcal{B}(M_1 \oplus \dots \oplus M_{\theta})$.
The Nichols algebra $\mathcal{B}(V)$ is known to be a~$\mathbb{N}_0^{\theta}$-graded algebra and coalgebra in ${}^H_H\mathcal{YD}$ such that $\deg
M_i=\alpha_i$ for all $i\in I$.

The adjoint action in the braided category ${}^H_H\mathcal{YD}$ is given by $\operatorname{ad}_cx(y)=xy-(x_{(-1)}\cdot y)x_{(0)}$ for all $x\in M_1
\oplus \cdots \oplus M_{\theta}$, $y\in \mathcal{B}(M)$, where $\delta(x)=x_{(-1)}\otimes x_{(0)}$.

By~\cite[Definition~6.8]{a-HeckSchn12a}, the Nichols algebra $\mathcal{B}(M)$ is called {\textit{decomposable}} if there exists
a~totally ordered index set $(L,\le)$ and a~family $(W_l)_{l\in L}$ of f\/inite-dimensional irreducible
$\mathbb{N}_0^\theta $-graded objects in $ {}^H_H\mathcal{YD}$ such that
\begin{gather}
\mathcal{B}(M)\simeq \bigotimes_{l\in L}\mathcal{B}(W_l).
\label{eq-decom}
\end{gather}
Decomposability of $\mathcal{B}(M)$ is known under several assumptions on~$M$.
In particular, if~$H$ is a~group algebra of an abelian group and $\dim M_i=1$ for all $1\leq i\leq \theta$, then
$\mathcal{B}(M)$ is decomposable by a~theorem of V.~Kharchenko~\cite[Theorem~2]{a-Khar99}.

Assume that $\mathcal{B}(M)$ is decomposable.
One def\/ines for any decomposition~\eqref{eq-decom} the set of positive roots $\boldsymbol{\Delta}^{[M]}_{+}\subset \mathbb{Z}^I$ and the set
of roots $\boldsymbol{\Delta}^{[M]}\subset \mathbb{Z}^I$ of $[M]$~by
\begin{gather*}
\boldsymbol{\Delta}^{[M]}_{+}=\{\deg(W_l)\,|\, l\in L\},
\qquad
\boldsymbol{\Delta}^{[M]}=\boldsymbol{\Delta}^{[M]}_{+}\cup-\boldsymbol{\Delta}^{[M]}_{+}.
\end{gather*}
The set of roots of $[M]$ does not depend on the choice of the decomposition~\eqref{eq-decom}.

Let $i\in I$.
Following~\cite[Definition~6.4]{HS10} we say that~$M$ is {\textit{$i$-finite}}, if for any $j\in I\setminus \{i\}$,
$(\ad_{c} M_i)^m (M_j)=0$ for some $m\in \mathbb{N}$.
Assume that~$M$ is~$i$-f\/inite.
Let $(a_{ij}^{M})_{j\in I}\in \mathbb{Z}^I$ and $R_i(M)=({R_i(M)}_j)_{j\in I}$, where
\begin{gather}
a_{ij}^M=
\begin{cases}
2& \text{if}\quad j=i,
\\
-\max \big\{m\in \mathbb{N}_0 \,|\, (\ad_c M_i)^m(M_j)\not=0 \big\}& \text{if}\quad j\not=i,
\end{cases}
\nonumber
\\
\label{eq-ri}
{R_i(M)}_i= {M_i}^*,
\qquad
{R_i(M)}_j= (\ad_{c}M_i)^{-a_{ij}^M}(M_j).
\end{gather}
Then $R_i(M)_j$ is irreducible by~\cite[Theorem~7.2(3)]{HS10}.
If~$M$ is not~$i$-f\/inite, then let $R_i(M)=M$.
Let
\begin{gather*}
\mathcal{X}_\theta^H(M)=\big\{[R_{i_1} \cdots R_{i_n}(M)]\in \mathcal{X}_\theta^H \,|\, n\in \mathbb{N}_0, \, i_1,\dots, i_n\in I\big\},
\\
\mathcal{F}_\theta^H(M)=\big\{R_{i_1} \cdots R_{i_n}(M)\in \mathcal{F}_\theta^H\,|\, n\in \mathbb{N}_0, \, i_1,\dots, i_n\in I\big\}.
\end{gather*}
We say that $M\in\mathcal{F}_\theta^H$ {\textit{admits all reflections}} if~$N$ is~$i$-f\/inite for all $N\in \mathcal{F}_\theta^H(M)$.

\begin{Theorem}
\label{theo.regualrcar}
Let $M\in \mathcal{F}_\theta^H$.
Assume that~$M$ admits all reflections.
Let $r: I\times \mathcal{X}_\theta^H(M)\rightarrow \mathcal{X}_\theta^H(M)$, $(i, [N])\mapsto [R_i(N)]$ for all $i\in I$.
Then
\begin{gather*}
\mathcal{C}(M)=\big\{I, \mathcal{X}_\theta^H(M), r, \big(A^{[N]}\big)_{[N]\in \mathcal{X}_\theta^H(M)}\big\}
\end{gather*}
is a~semi-Cartan graph.
If moreover $\rersys{[M]}$ is finite, then $\mathcal{C}(M)$ is a~Cartan graph.
\end{Theorem}
\begin{proof}
The f\/irst claim follows from~\cite[Theorem~3.12]{a-AHS08}, see~\cite[Theorem~6.10]{HS10} for details.
If~$\rersys{[M]}$ is f\/inite and~$M$ admits all ref\/lections, then
\begin{gather*}
\mathcal{R}(M)=\mathcal{R}\big(\mathcal{C}(M),\big(\rersys{X}\big)_{X\in \mathcal{X}_\theta^H(M)}\big)
\end{gather*}
is a~root system of type $\mathcal{C}(M)$ by~\cite[Corollary~6.16]{a-HeckSchn12a}.
Hence $\mathcal{C}(M)$ is a~Cartan graph because of~\cite[Proposition~2.9]{c-Heck09b}.
\end{proof}

Therefore if~$M$ admits all ref\/lections then we can attach the groupoid $\mathcal{W}(M):=\mathcal{W}(\mathcal{C}(M))$ to~$M$.

\subsection{Small Cartan graphs for Nichols algebras of diagonal type}
\label{se:dynkin diag}
Let~$G$ be an abelian group and let~$V$ be a~Yetter--Drinfel'd module over $\Bbbk G$ of rank~$\theta$ with a~basis $\{x_i| i\in I \}$.
Let $.: \Bbbk G\otimes V\rightarrow V$ and $\delta: V\rightarrow \Bbbk G\otimes V$ denote the left action and the left
coaction of $\Bbbk G$ on~$V$, respectively.
Assume that~$V$ is of diagonal type.
More precisely, let $\{g_i \,|\, i\in I\}$ be a~subset of~$G$ and $q_{ij}\in \Bbbk^\ast$ for all $i,j\in I$, such that
\begin{gather*}
\delta (x_i)=g_i\otimes x_i,
\qquad
g_i\lact x_j=q_{ij}x_j
\end{gather*}
for $i,j\in I$.
Then~$V$ is a~{\textit{braided vector space}} of dimension~$\theta$~\cite[Definition~5.4]{inp-Andr02} and the braiding
$c \in \mathrm{End}(V\otimes V)$ is of diagonal type, that is $c(x_i\otimes x_j)=q_{ij}x_j\otimes x_i$ for all $i,j\in I$.
Then $(q_{ij})_{i,j\in I}$ is the {\textit{braiding matrix of~$V$}} with respect to the basis $\{x_i| i\in I \}$.

The braiding matrix is known to be independent of the basis $\{x_i|i\in I\}$, up to permutation of~$I$.
It can be obtained for example from~\cite[Proposition~1.3]{HS101} using the arguments in its proof.

The Nichols algebra $\mathcal{B}(V)$ generated by~$V$ is said to be {\textit{of diagonal
type}}~\cite[Definition~5.8]{inp-Andr02}.

For $\rho \in \mathcal{B}(V)$, the braided commutator $\ad_{c}$ takes the form $\ad_{c}x_i(\rho)=x_i\rho-(g_i.\rho) x_i$.
\begin{Lemma}
\label{le:aijM}
Let $M=(M_i)_{i\in I}\in \mathcal{F}_\theta^{\Bbbk G}$ be a~tuple of one-dimensional Yetter--Drinfel'd modules over $\Bbbk G$ and let $x_i$ be a~basis of
$M_i$, for all $i\in I$.
Let $m\in \mathbb{N}_0$.
For any $i,j\in I$ with $i\not=j$, the following are equivalent:
\begin{enumerate}\itemsep=0pt
\item[$(a)$] $(m+1)_{q_{ii}}(q_{ii}^mq_{ij}q_{ji}-1)=0$ and $(k+1)_{q_{ii}}(q_{ii}^kq_{ij}q_{ji}-1)\not=0$ for all $0\leq k<m$,
\item[$(b)$] $(\ad_{c}x_i)^{m+1}(x_j)=0$ and $(\ad_{c}x_i)^m(x_j)\not=0$ in $\mathcal{B}(V)$,
\item[$(c)$] $-a_{ij}^M=m$.
\end{enumerate}
Here $(n)_q:=1+q+\dots+q^{n-1}$, which is $0$ if and only if $q^n=0$ for $q\not=1$ or $p|n$ for $q=1$.
\end{Lemma}
\begin{proof}
(a)$\Leftrightarrow$(b) follows from~\cite[Lemma~3.7]{inp-AndrSchn02} and (b)$\Leftrightarrow$(c) holds by the
def\/inition of $a_{ij}^M$.
\end{proof}
\begin{Lemma}
\label{le:ifinite}
Let $i\in I$.
Then $M=(M_j)_{j\in I}\in \mathcal{F}_\theta^{\Bbbk G}$ is~$i$-finite if and only if for any $j\in I\setminus\{i\}$ there is a~non-negative
integer~$m$ satisfying $(m+1)_{q_{ii}}(q_{ii}^mq_{ij}q_{ji}-1)=0$.
\end{Lemma}
\begin{proof}
The claim follows from Lemma~\ref{le:aijM}.
\end{proof}

Let~$V$ be a~$\theta$-dimensional braided vector space of diagonal type.
Let $(q_{ij})_{i,j\in I}$ be a~braiding matrix of~$V$.
The \textit{Dynkin diagram}~\cite[Definition~4]{a-Heck04e} of~$V$ is denoted by $\mathcal{D}$.
It is a~non-directed graph with the following properties:
\begin{itemize}\itemsep=0pt
\item there is a~bijective map~$\phi$ from~$I$ to the vertices of $\mathcal{D}$,
\item for all $i\in I$ the vertex $\phi (i)$ is labeled by $q_{ii}$,
\item for any $i,j\in I$ with $i\not=j$, the number $n_{ij}$ of edges between $\phi (i)$ and $\phi (j)$ is either $0$
or~$1$.
If $q_{ij}q_{ji}=1$ then $n_{ij}=0$, otherwise $n_{ij}=1$ and the edge is labeled by $q_{ij}q_{ji}$.
\end{itemize}

Let $M=(M_i)_{i\in I}\in \mathcal{F}_\theta^{\Bbbk G}$ be a~tuple of one-dimensional Yetter--Drinfel'd modules over $\Bbbk G$.
The \textit{Dynkin diagram of~$M$} is the Dynkin diagram of the braided vector space $M_1\oplus\dots \oplus M_{\theta}$.

Let $i\in I$.
Assume that~$M$ is~$i$-f\/inite.
By def\/initions of $R_i(M)$ and~$M$, the tuple $(y_{j})_{j\in I}$ is a~basis of $R_i(M)$, where
\begin{gather*}
y_j:=
\begin{cases}
y_i & \text{if}\quad j=i,
\\
(\ad_{c} x_i)^{-a_{ij}^M}(x_j) & \text{if}\quad j\not=i,
\end{cases}
\end{gather*}
where $y_i\in M_i^*\setminus \{0\}$.

From the method in~\cite[Example 1]{a-Heck04e}, one can obtain the labels of the Dynkin diagram of
$R_i(M)=(R_i(M)_j)_{j\in I}$.
In more detail, we have the following lemma.

\begin{Lemma}
\label{le:Dynkin}
Let $i\in I$.
Assume that~$M$ is~$i$-finite and let $a_{ij}:=a_{ij}^M$ for all $j\in I$.
Let $(q'_{jk})_{j,k\in I}$ be the braiding matrix of $R_i(M)$ with respect to $(y_j)_{j\in I}$.
Then the labels of the Dynkin diagram of $R_i(M)=(R_i(M)_j)_{j\in I}$ are
\begin{gather*}
q_{jj}'=
\begin{cases}
q_{ii} & \text{if}\quad j=i,
\\
q_{jj} & \text{if}\quad j\not=i,\quad q_{ij}q_{ji}=q_{ii}^{a_{ij}},
\\
q_{ii}q_{jj}{(q_{ij}q_{ji})^{-a_{ij}}} & \text{if}\quad j\not=i,\quad q_{ii}\in G'_{1-a_{ij}},
\\
q_{jj}{(q_{ij}q_{ji})^{-a_{ij}}} & \text{if}\quad j\not=i,\quad q_{ii}=1,
\end{cases}
\\
q_{ij}'q_{ji}'=
\begin{cases}
q_{ij}q_{ji} & \text{if}\quad j\not=i,\quad q_{ij}q_{ji}=q_{ii}^{a_{ij}},
\\
q_{ii}^2(q_{ij}q_{ji})^{-1} & \text{if}\quad j\not=i,\quad q_{ii}\in G'_{1-a_{ij}},
\\
(q_{ij}q_{ji})^{-1} & \text{if}\quad j\not=i,\quad q_{ii}=1,
\end{cases}
\end{gather*}
and if $j, k\not=i$, $j\not=k$, then
\begin{gather*}
q_{jk}'q_{kj}'=
\begin{cases}
q_{jk}q_{kj} & \text{if}\quad q_{ir}q_{ri}=q_{ii}^{a_{ir}},\quad r\in \{j, k\},
\\
q_{jk}q_{kj}(q_{ik}q_{ki}q_{ii}^{-1})^{-a_{ij}}& \text{if}\quad q_{ij}q_{ji}=q_{ii}^{a_{ij}},\quad q_{ii}\in G'_{1-a_{ik}},
\\
q_{jk}q_{kj}(q_{ij}q_{ji})^{-a_{ik}}(q_{ik}q_{ki})^{-a_{ij}} & \text{if}\quad q_{ii}=1,
\\
q_{jk}q_{kj}q_{ii}^{2}(q_{ij}q_{ji}q_{ik}q_{ki})^{-a_{ij}} & \text{if}\quad q_{ii}\in G'_{1-a_{ik}},\quad q_{ii}\in G'_{1-a_{ij}}.
\end{cases}
\end{gather*}
\end{Lemma}
Assume that~$M$ admits all ref\/lections.
By Theorem~\ref{theo.regualrcar}, we are able to construct a~semi-Cartan graph $\mathcal{C}(M)$ of~$M$
\begin{gather*}
\mathcal{C}(M)=\big(I,\mathcal{X}_\theta^{\Bbbk G}(M), (r_i)_{i\in I}, \big(A^X\big)_{X\in \mathcal{X}_\theta^{\Bbbk G}(M)}\big),
\end{gather*}
where $\mathcal{X}_\theta^{\Bbbk G}(M)=\{[R_{i_1}\cdots R_{i_n}(M)]\in \mathcal{F}_\theta^{\Bbbk G}\,|\, n\in \mathbb{N}_0,\, i_1, \dots, i_n\in I\}$.
Note that any $X\in \mathcal{X}_\theta^{\Bbbk G}(M)$ has a~well-def\/ined braiding matrix given by the braiding matrix of any representative of~$X$.
\begin{Definition}
Assume that~$M$ admits all ref\/lections.
For all $X\in \mathcal{X}_\theta^{\Bbbk G}(M)$ let
\begin{gather*}
[X]_{\theta}^s=\big\{Y\in \mathcal{X}_\theta^{\Bbbk G}(M)\,\big|\, \text{$Y$ and $X$ have the same Dynkin diagram}\big\}.
\end{gather*}
Let $\mathcal{Y}_{\theta}^s(M)=\{[X]^s_{\theta} | X\in \mathcal{X}_\theta^{\Bbbk G}(M)\}$ and $A^{[X]^s_{\theta}}=A^X$ for all $X\in
\mathcal{X}_\theta^{\Bbbk G}(M)$.
Let $t: I\times \mathcal{Y}_{\theta}^s(M)\rightarrow \mathcal{Y}_{\theta}^s(M)$, $(i, [X]^s_{\theta})\mapsto [r_i(X)]^s_{\theta}$.
Then the tuple
\begin{gather*}
\mathcal{C}_s(M)=\big\{I, \mathcal{Y}_{\theta}^s(M), t, \big(A^Y\big)_{Y\in \mathcal{Y}_{\theta}^s(M)}\big\}
\end{gather*}
is called the small semi-Cartan graph of~$M$.
\end{Definition}

\begin{Proposition}
\label{prop:Cs}
Assume that~$M$ admits all reflections.
Then the tuple
\begin{gather*}
\mathcal{C}_s(M)=\big\{I, \mathcal{Y}_{\theta}^s(M), t, \big(A^Y\big)_{Y\in \mathcal{Y}_{\theta}^s(M)}\big\}
\end{gather*}
is a~semi-Cartan graph.
Moreover, if $\mathcal{C}(M)$ is a~finite Cartan graph, then $\mathcal{C}_s(M)$ is a~finite Cartan graph.
\end{Proposition}
\begin{proof}
The map~$t$ and $A^{[X]^s_{\theta}}$ are well-def\/ined for all $[X]^s_{\theta}\in \mathcal{Y}_{\theta}^s(M)$.
Indeed, if $X, X' \in [X]^s_{\theta} $, then $A^X=A^{X'}$ by Lemma~\ref{le:aijM}.
Thus $A^{[X]^s_{\theta}}$ is well-def\/ined.
Further, Lemma~\ref{le:Dynkin} implies that $r_i(X)$ and $r_i(X')$ have the same Dynkin diagram.
Hence~$t$ is well-def\/ined.
Since $t_i([X]^s_{\theta})=t(i, [X]^s_{\theta})=[r_i(X)]^s_{\theta}$ and $r_i^2=\mathrm{id}$, then for all $i\in I$,
$t_i^2=\mathrm{id}_{\mathcal{Y}_{\theta}^s(M)}$.
Moreover, $a_{ij}^{[X]^s_{\theta}}=a_{ij}^{t_i([X]^s_{\theta})}$. Hence the f\/irst claim holds.

The second claim follows from the def\/initions of $\rersys{X}$ and $\mathcal{C}_s(M)$.
\end{proof}

\section{Finite Cartan graphs of rank two}
\label{se:classi of Carscheme}
In this section we simplify slightly the fundaments of the theory presented in~\cite{c-Heck09a} and give
a~characterization of f\/inite Cartan graphs of rank two.
\begin{Definition}
Let $\mathcal{A}^+$ denote the smallest subset of $\cup_{n\geq2}\mathbb{N}_0^n$ such that
\begin{itemize}\itemsep=0pt
\item $(0, 0)\in \mathcal{A}^+$,
\item if $(c_1, \dots, c_n)\in \mathcal{A}^+$ and $1<i\leq n$, then $(c_1, \dots, c_{i-2}, c_{i-1}+1, 1, c_i+1, \dots, c_n)\in
\mathcal{A}^+$.
\end{itemize}
\end{Definition}
\begin{Remark}
Note that our def\/inition of $\mathcal{A}^+$ is dif\/ferent from the one in~\cite{c-Heck11a}.
\end{Remark}

From the def\/inition of $\mathcal{A}^+$, we get the following lemma.
\begin{Lemma}
\label{le:sumA}
Let $n\geq 2$ and $(c_1, \dots,c_n)\in \mathcal{A}^+$.
Then $\sum\limits_{i=1}^{n}c_i=3n-6$.
\end{Lemma}
The def\/inition of $\mathcal{A}^+$ implies the following.

\begin{Proposition}
\label{prop:triangula}
Let $n\geq 2$.
Enumerate the vertices of a~convex~$n$-gon by $1, \dots, n$ such that consecutive integers correspond to neighboring
vertices.
Let $T_n$ be the set of triangulations of a~convex~$n$-gon with non-intersecting diagonals.
Let $T=\cup_{n\geq 2}T_n$.
For any triangulation $t\in T_n$ and any $i\in \{1, \dots, n\}$, let $c_i$ be the number of triangles meeting at
the~$i$-th vertex.
Then the map $\psi: T \rightarrow \mathcal{A}^+$, $t\mapsto (c_1, \dots, c_n)$ is a~bijection.
\end{Proposition}
\begin{proof}
We proceed by induction on~$n$.
For $n=2$, a~triangulation of a~convex $2$-gon is itself.
Then $(c_1, c_2)=(0, 0)$.
Hence the claim is true for $n=2$.
For $n\geq 3$, the def\/inition of $\mathcal{A}^+$ corresponds bijectively to the construction of a~triangulation of a~convex
$(n+1)$-gon by adding a~new triangle between two consecutive vertices of a~convex~$n$-gon, but not at the edge between
the f\/irst and the last vertex.
By adding one triangle between two consecutive vertices of a~convex~$n$-gon, one increases the number of triangles at
the two adjacent vertices and the number of triangles at the new vertex becomes~$1$.
\end{proof}

\begin{Corollary}\samepage
\label{coro:A+}
Let $n\geq 2$ and let $(c_1, \dots,c_n)\in \mathcal{A}^+$.
\begin{enumerate}\itemsep=0pt
\item[$(1)$] $(c_n, c_{n-1}, \dots,c_1)\in \mathcal{A}^+$ and $(c_2, c_3, \dots, c_n, c_1)\in \mathcal{A}^+$.
\item[$(2)$] If $n\geq3$, then there is $1<i<n$ satisfying $c_i=1$.
For any such~$i$, $(c_1, \dots, c_{i-2}, c_{i-1}-1,$ $c_{i+1}-1, c_{i+2}, \dots, c_n)\in \mathcal{A}^+$.
\item[$(3)$] If $n\geq 3$, then $c_i\geq 1$ for all $1\leq i\leq n$.
\item[$(4)$] If $c_i=1$, $c_{i+1}=1$ for some $1\leq i\leq n-1$, then $n=3$ and $c=(1, 1, 1)$.
\end{enumerate}
\end{Corollary}

\begin{proof}
$(1)$ and $(2)$ follow directly from the bijection between $\mathcal{A}^+$ and triangulations of convex~$n$-gons in
Proposition~\ref{prop:triangula}.
$(3)$ follows from the def\/inition of $\mathcal{A}^+$.
$(4)$ follows from $(2)$ and $(3)$.
\end{proof}
We say that two consecutive entries of a~sequence in $\mathcal{A}^+$ are neighbors.
\begin{Theorem}
\label{Theo:clascar}
Let $n\geq3$.
Then any sequence $(c_1, \dots,c_n)\in \mathcal{A}^+$ contains a~subsequence $(c_k)_{i \leq k\leq j}$, where $1\leq i \leq j
\leq n$, of the form
\begin{gather*}
(1,1),
\qquad
(1, 2, a),
\qquad
(2, 1, b),
\qquad
(1, 3, 1, b)
\end{gather*}
or their transpose, where $1\leq a\leq 3$ and $3\leq b\leq 5$.
\end{Theorem}
\begin{Remark}
We record that it is natural to exclude the cases $b=1$ and $b=2$ since $(1,3,1,1)$ contains the subsequence $(1,1)$ and
$(1,3,1,2)$ contains the transpose of $(2,1,3)$.
\end{Remark}
\begin{Remark}
The claim becomes false by omitting one of the sequences from the theorem.
In Table~\ref{table:star} we list sequences in $\mathcal{A}^+$ which contain precisely one of the sequences in
Theorem~\ref{Theo:clascar}.
\begin{table}[t] \centering
\caption{Sequences in $\mathcal{A}^+$ containing exactly one subsequence from Theorem~\ref{Theo:clascar}.}\label{table:star}
\vspace{1mm}
\begin{tabular}{r|r}
\hline
\text{subsequences} & sequences in $\mathcal{A}^+$
\\
\hline
\hline
$(1,1)$ & $(1,1,1)$
\\
\hline
$(1,2,1)$ &$(1,2,1,2)$
\\
\hline
$(1,2,2)$ &$(1,2,2,2,2,2,1,6)$
\\
\hline
$(1,2,3)$ &$(1,2,3,1,6,1,2,3,1,6,1,2,3,1,6)$
\\
\hline
$(2,1,3)$ &$(2,1,3,4,2,1,3,4,2,1,3,4)$
\\
\hline
$(2,1,4)$ &$(2,1,4,2,1,4,2,1,4)$
\\
\hline
$(2,1,5)$ &$(2,1,5,1,2,4,2,1,5,1,2,4)$
\\
\hline
$(1,3,1,3)$ & $(1,3,1,3,1,3)$
\\
\hline
$(1,3,1,4)$ & $(1,3,1,4,1,3,1,4)$
\\
\hline
$(1,3,1,5)$ & $(1,3,1,5,1,3,1,5,1,3,1,5)$
\\
\hline
\end{tabular}
\end{table}
\end{Remark}

\begin{proof}
Let $c=(c_1, \dots,c_n)\in \mathcal{A}^+$ such that the claim does not hold for~$c$.
Then $n\geq5$ and~$c$ has no subsequence $(2,1,2)$.
Otherwise $c=(1, 2, 1, 2)$ or $c=(2, 1, 2, 1)$ by Corollary~\ref{coro:A+}(2),(4).
We def\/ine $E=\{\nu_{ij}\,|\, i,j\in \{1,2\}\}$, where the sequences $\nu_{ij}$ are given by
\begin{gather*}
\nu_{11}=(1),
\qquad
\nu_{12}=(2,1),
\qquad
\nu_{21}=(1,2),
\qquad
\nu_{22}=(1,3,1).
\end{gather*}

Now we decompose~$c$ by the following steps.

Replace all subsequences $(2,1)$ by $\nu_{12}$, then all subsequences $(1,2)$ by $\nu_{21}$, then all subsequences
$(1,3,1)$ by $\nu_{22}$, and f\/inally all entries $1$ by $\nu_{11}$.
By this construction, $(\nu_{11}, 3, \nu_{11})$ is not a~subsequence of~$d$.
Hence we get a~decomposition $d=(d_1,\dots,d_k)$, where $k\geq2$, of~$c$ into subsequences of the form $(a)$ and~$\nu$,
where $a\geq2$ and $\nu\in E$.

Since $(1,1)$, $(1, 2, a)$, $(2,1,b)$, $(1,3,1,b)$ and their transposes are not subsequences of~$c$, where $1\leq a~\leq
3$, $2\leq b\leq 5$, we obtain the following conditions on the entries of~$d$:
\begin{itemize}\itemsep=0pt
\item[--] no entry $\nu_{ij}$ of~$d$, where $i,j\in \{1,2\}$, has $2$ or $\nu_{kl}$ with $k,l\in \{1,2\}$ as a~neighbor.
\item[--] if $(\nu_{21}, a)$ or $(a, \nu_{12})$ is a~subsequence of~$d$, then $a\geq4$.
\item[--] if $(\nu_{i2}, b)$ or $(b, \nu_{2i})$ is a~subsequence of~$d$, where $i\in \{1,2\}$, then $b\geq6$.
\end{itemize}

By applying Corollary~\ref{coro:A+}(2) we get further reductions of~$d$:
\begin{gather*}
(\dots, d_{m-1}, \nu_{ij}, d_{m+1}, \dots)\rightarrow (\dots, d_{m-1}\!-i, d_{m+1}\!-j),
\qquad\!
(\nu_{i2}, d_2, \dots)\rightarrow (\nu_{i1}, d_2-1, \dots),
\end{gather*}
where $i,j\in \{1,2\}$.
Thus we can perform such reductions at all places in~$d$, where an entry $\nu_{ij}$ with $i,j\in\{1,2\}$ appears.
After decreasing them, we get $d_m\geq2$, where $1<m< k$.
Indeed, we get the following conditions.
\begin{itemize}\itemsep=0pt
\item[--] If $d=(\dots, d_{m-1}, d_m, d_{m+1}, \dots)$, where $d_m\geq 6$, then $d_m$ can be reduced at most by~$4$.
Hence the value of $d_m$ after reduction is at least~$2$.
\item[--] If $4\leq d_m\leq 5$, then neither $(\nu_{i2}, d_m)$ nor $(d_m, \nu_{2i})$ is a~subsequence of~$d$, where
$i\in \{1,2\}$.
Hence~$d$ can be reduced by at most~$2$.
\item[--] If $d_m=3$, then $d_{m-1}, d_{m+1}\notin \{\nu_{12}, \nu_{21}, \nu_{22}\}$.
Further, $(d_{m-1}, d_{m+1})\not=(\nu_{11}, \nu_{11})$.
Hence $d_m$ decreases by at most~$1$.
\item[--] If $d_m=2$, then it has no neighbour $\nu_{ij}$ with $i, j\in \{1,2\}$.
Hence $d_m$ does not change.
\end{itemize}
Thus one can reduce~$c$ to a~sequence $(c_1', \dots, c_l')$ with $l\geq1$, where $c_m'\geq 2$ for all $1<m<l$ and $c_1',
c_l'\geq 1$.
This is a~contradiction to Corollary~\ref{coro:A+}(2).
\end{proof}
Recall that $(\alpha_1, \alpha_2)$ is the standard basis of $\mathbb{Z}^2$.
We def\/ine a~map
\begin{gather}
\label{eq-etadefi}
\eta: \ \mathbb{Z}\rightarrow \SLZ{2},
\qquad
a\mapsto
\begin{pmatrix}
a & -1
\\
1 & 0
\end{pmatrix}.
\end{gather}

\begin{Lemma}
\label{le:beta}
Let $n\in \mathbb{N}$ and $(c_k)_{1\leq k\leq n}\in \mathbb{Z}^n$.
For all $1\leq k\leq n+1$, let $\beta_0=-\alpha_2$ and $\beta_k=\eta(c_1)\dots\eta(c_{k-1})(\alpha_1)$.
Then the following hold:
\begin{enumerate}\itemsep=0pt
\item[$(1)$] $\beta_{k+1}=c_k\beta_k-\beta_{k-1}$ for all $1\leq k\leq n$,
\item[$(2)$] if $c_1\geq 1$ and $c_k\geq 2$ for all $1<k<n$, then $\beta_k\in \mathbb{N}_0^2$ for all $1\leq k\leq n$ and
$\beta_{k}-\beta_{k-1}\in \mathbb{N}_0^2\setminus \{0\}$ for $1<k\leq n$.
\end{enumerate}
\end{Lemma}
\begin{proof}
(1) By def\/inition, $\beta_1=\alpha_1$ and $\beta_2=\eta(c_1)(\alpha_1)=c_1\alpha_1+\alpha_2$.
Thus the claim holds for $k=1$.
Since $\eta(c_{k-1})(\alpha_2)=-\alpha_1$, then
\begin{gather*}
\beta_{k+1}= \eta(c_1)\dots \eta(c_k)(\alpha_1)
= \eta(c_1)\dots \eta(c_{k-1})(c_k\alpha_1+\alpha_2)
= c_k\beta_k-\beta_{k-1}
\end{gather*}
for all $k\geq2$.

(2) For all $0\leq k \leq n$, let $a_k$, $b_k\in \mathbb{Z}$ such that $\beta_k=a_k\alpha_1+b_k\alpha_2$.
By induction on~$k$, we get the following.
\begin{itemize}\itemsep=0pt
\item If $c_k\geq 2$ for $1\leq k< n$, then
\begin{gather*}
a_k> b_k\geq 0,
\qquad
a_k>a_{k-1},
\qquad
b_k>b_{k-1},
\qquad
a_k-b_k-(a_{k-1}-b_{k-1})\geq 0
\end{gather*}
for all $1\leq k\leq n$.
\item If $c_1=1$ and $c_k\geq 2$ for $2\leq k< n$, then
\begin{gather*}
b_k\geq a_k> 0,
\qquad
a_k\geq a_{k-1},
\qquad
b_k>b_{k-1},
\qquad
a_k-b_k-(a_{k-1}-b_{k-1})<0.
\end{gather*}
for all $2\leq k\leq n$.
\end{itemize}
Thus $\beta_k\in \mathbb{N}_0^2$ for all $1\leq k\leq n$ and $\beta_k-\beta_{k-1}\in \mathbb{N}_0^2\setminus \{0\}$ for $1<k\leq n$.
\end{proof}
The following theorem will be used in the proof of Theorem~\ref{theo.cartan}.
It was proven partially in~\cite[Propositon~5.3]{c-Heck09a}.
Notice that the def\/inition of $\mathcal{A}^+$ is dif\/ferent from the one in~\cite{c-Heck09a}.

{\samepage
\begin{Theorem}
\label{theo:cqid}
Let $n\geq2$ and $(c_i)_{1\leq i\leq n}\in \mathbb{Z}^n$.
Then the following are equivalent:
\begin{itemize}\itemsep=0pt
\item[$(1)$] $(c_i)_{1\leq i\leq n}\in \mathcal{A}^+$,
\item[$(2)$] $\eta(c_1)\cdots\eta(c_n)=-\mathrm{id}$ and $\beta_k=\eta(c_1)\cdots\eta(c_{k-1})(\alpha_1)\in \mathbb{N}_0^2$ for
all $1\leq k\leq n$.
\end{itemize}
\end{Theorem}

}

\begin{proof}
(1)$\Rightarrow$(2).
We apply induction on~$n$.
If $n=2$, then $(c_1, c_2)=(0, 0)$, $\eta(0)^2=-\mathrm{id}$, $\beta_1=\alpha_1$, and $\beta_2=\alpha_2$.
Assume that $n\geq3$.
By the def\/inition of $\mathcal{A}^+$, there is a~$(c_1', \dots, c_{n-1}')\in \mathcal{A}^+$ and $1<i\leq n-1$ such that
\begin{gather*}
(c_1, \dots, c_n)=(c_1',\dots, c_{i-1}'+1, 1,c_{i}'+1,c_{i+1}',\dots,c_{n-1}').
\end{gather*}
By calculation,
\begin{gather}
\label{eq-eta}
\eta(a)\eta(b)=\eta(a+1)\eta(1)\eta(b+1)
\end{gather}
for all $a, b\in \mathbb{Z}$.
Then
\begin{gather*}
\eta(c_1)\cdots\eta(c_n)=\eta(c_1')\cdots\eta(c_{n-1}')=-\mathrm{id}.
\end{gather*}
Let $\beta_i'=\eta(c_1')\cdots\eta(c_{i-1}')(\alpha_1)$ for all $1\leq i \leq n-1$.
Then $\beta_k=\beta_k'$ for all $1\leq k<i$ and $\beta_{k}=\beta_{k-1}'$ for all $i+1 \leq k\leq n+1$.
Finally
\begin{gather*}
\beta_i= \eta(c_1)\cdots\eta(c_{i-1})(\alpha_1)
= \eta(c_1')\cdots\eta(c_{i-2}')\eta(c_{i-1}'+1)(\alpha_1)
\\
\phantom{\beta_i}
= \eta(c_1')\cdots\eta(c_{i-2}')(\eta(c_{i-1}')(\alpha_1)+\alpha_1)
= \beta_i'+\beta_{i-1}'\in \mathbb{N}_0^2.
\end{gather*}
Then (2) follows.

(2)$\Rightarrow$(1).
Again we proceed by induction.
If $n=2$, then
\begin{gather*}
\eta(c_1)\eta(c_2)=
\begin{pmatrix}
c_1c_2-1& -c_1
\\
c_2&-1
\end{pmatrix}
=-\mathrm{id}
\end{gather*}
implies that $(c_1, c_2)=(0, 0)\in \mathcal{A}^+$.
Assume that $n\geq3$.
Set $\beta_0=-\alpha_2$.
One has $\beta_{k+1}=c_k\beta_k-\beta_{k-1}$ for all $1\leq k<n$.
By assumption, the condition $\beta_{k-1}, \beta_k, \beta_{k+1}\in \mathbb{N}_0^2$ implies $c_k>0$ for $2\leq k<n$ and
$c_1\geq0$.
If $c_1=0$ then $\beta_2=\alpha_2$ and $\beta_3= c_2\alpha_2-\alpha_1\notin \mathbb{N}_0^2$.
Hence $c_k\geq 1$ for all $1\leq k<n$.
Moreover, there is $1<i<n$ satisfying $c_i=1$.
Indeed, $\beta_{n+1}=c_n\beta_{n}-\beta_{n-1}=(c_n-1)\beta_n+(\beta_n-\beta_{n-1})$ by Lemma~\ref{le:beta}(1).
Assume that $c_i\geq 2$ for all $1<i<n$.
Then $\beta_{n+1}\in \mathbb{N}_0^2$ if $c_n\geq 1$ and $-\beta_{n+1}\in \mathbb{N}_0^2\setminus \{0, \alpha_1\}$ if $c_n\leq 0$~by
Lemma~\ref{le:beta}(2), since $n\geq 3$.
This is a~contradiction to $\beta_{n+1}=\eta(c_1)\cdots\eta(c_n)(\alpha_1)=(-\mathrm{id})(\alpha_1)=-\alpha_1$.

Hence there is $(c_1', \dots, c_{n-1}')\in \mathbb{Z}^{n-1}$ such that
\begin{gather*}
(c_1, \dots, c_n)=(c_1', \dots, c_{i-1}'+1, 1, c_i'+1, c_{i+1}', \dots, c_{n-1}').
\end{gather*}
Then $\eta(c_1)\cdots\eta(c_n)=\eta(c_1')\cdots\eta(c_{n-1}')=-\mathrm{id}$ by equation~\eqref{eq-eta} and
$\beta_k'=\eta(c_1')\dots\eta(c_{k-1}')(\alpha_1)\in \mathbb{N}_0^2$ for all $1\leq k\leq n-1$.
Hence $(c_1', \dots, c_{n-1}')\in \mathcal{A}^+$ by induction hypothesis.
Then
\begin{gather*}
(c_1, \dots, c_n)\in \mathcal{A}^+.
\tag*{\qed}
\end{gather*}
\renewcommand{\qed}{}
\end{proof}
\begin{Definition}
Let $\mathcal{C}= \mathcal{C}(I, \mathcal{X}, r, (A^X)_{X\in \mathcal{X}})$ be a~semi-Cartan graph of rank two and let $X\in \mathcal{X}$ and $i\in I$.
The \textit{characteristic sequence of $\mathcal{C}$ with respect to~$X$ and~$i$} is the inf\/inite sequence $(c_{k}^{X,i})_{k\geq
1}$ of non-negative integers, where
\begin{gather*}
c_{2k+1}^{X,i}=-a_{ij}^{(r_jr_i)^k(X)}=-a_{ij}^{r_i(r_jr_i)^k(X)},
\qquad
c_{2k+2}^{X,i}=-a_{ji}^{r_i(r_jr_i)^k(X)}=-a_{ji}^{(r_jr_i)^{k+1}(X)}
\end{gather*}
for all $k\geq 0$ and $j\in I\backslash \{i\}$.
\end{Definition}
By the def\/inition of a~characteristic sequence, we get the following remark.
\begin{Remark}
\label{re:char_se_a_j}
Let $\mathcal{C}= \mathcal{C}(I, \mathcal{X}, r, (A^X)_{X\in \mathcal{X}})$ be a~semi-Cartan graph of rank two and let $X\in \mathcal{X}$ and $i,j\in I$ with
$i\not=j$.
Let $(c_{k})_{k\geq 1}$ be the characteristic sequence of $\mathcal{C}$ with respect to~$X$ and~$i$.
\begin{itemize}\itemsep=0pt
\item[--] The characteristic sequence of $\mathcal{C}$ with respect to $r_i(X)$ and~$j$ is $(c_{k+1})_{k\geq1}$.
\item[--] Suppose that $(r_jr_i)^n(X)=X$ for some $n\geq 1$.
Then the characteristic sequence of $\mathcal{C}$ with respect to~$X$ and~$j$ is $(c_{2n+1-k})_{k\geq1}$.
\end{itemize}
\end{Remark}
\begin{Definition}
Let $\mathcal{C}= \mathcal{C}(I, \mathcal{X}, r, (A^X)_{X\in \mathcal{X}})$ be a~semi-Cartan graph of rank two and let $X\in \mathcal{X}$ and $i\in I$.
Let $(c_{k})_{k\geq 1}$ be the characteristic sequence of $\mathcal{C}$ with respect to~$X$ and~$i$.
The \textit{root sequence of $\mathcal{C}$ with respect to~$X$ and~$i$}is the inf\/inite sequence $(\beta_k)_{k\geq 1}$ of
elements of $\mathbb{Z}^2$, where
\begin{gather*}
\beta_k=\eta(c_1)\cdots\eta(c_{k-1})(\alpha_1)
\end{gather*}
for all $k\geq 1$.
In particular, $\beta_1=\alpha_1$.
\end{Definition}

Let $\mathcal{C}=\mathcal{C}(I=\{1, 2\}, \mathcal{X}, r, (A^X)_{X\in \mathcal{X}})$ be a~semi-Cartan graph.
For all $X\in \mathcal{X}$, the maps~$s_1^X$,~$s_2^X$ are def\/ined by equation~\eqref{eq-si}.
Recall that $(\alpha_1, \alpha_2)$ is a~basis of $\mathbb{Z}^2$ and~$\eta$ is a~map def\/ined by equation~\eqref{eq-etadefi}.
Def\/ine a~map
\begin{gather*}
\tau: \ \mathbb{Z}^2\rightarrow \mathbb{Z}^2,
\qquad
a\alpha_1+b\alpha_2 \mapsto b\alpha_1+a\alpha_2
\end{gather*}
for any $a, b\in \mathbb{Z}$.
One obtains
\begin{gather}
\label{eq-sitau}
s_1^X=\eta\big({-}a_{12}^X\big)\tau,
\qquad
s_2^X=\tau\eta\big({-}a_{21}^X\big)
\end{gather}
for all $X\in \mathcal{X}$.

\begin{Lemma}
\label{le:root}
Let $\mathcal{C}=\mathcal{C}(I=\{1, 2\}, \mathcal{X}, r, (A^X)_{X\in \mathcal{X}})$ be a~semi-Cartan graph of rank two and let $X\in \mathcal{X}$.
Let $(\beta_{k})_{k\geq 1}$ be the root sequence of $\mathcal{C}$ with respect to~$X$ and $1$ and let $(\gamma_{k})_{k\geq 1}$
be the root sequence of $\mathcal{C}$ with respect to~$X$ and~$2$.
Then
\begin{alignat*}{3}
& \beta_{2k+1}= \mathrm{id}_X(s_1s_2)^k \alpha_1,
\qquad &&
\beta_{2k+2}=\mathrm{id}_X(s_1s_2)^k s_1 \alpha_2,&
\\
& \tau \gamma_{2k+1}= \mathrm{id}_X(s_2 s_1)^k \alpha_2,
\qquad &&
\tau \gamma_{2k+2}=\mathrm{id}_X(s_2 s_1)^k s_2 \alpha_1 &
\end{alignat*}
for all $k\geq 0$.
Hence $\rersys{X}=\{\pm \beta_k, \pm\tau\gamma_k|k\geq 1\}$.
\end{Lemma}

\begin{proof}
Let $(c_k)_{k\geq 1}$ be the characteristic sequence with respect to~$X$ and $i=1$.
By equation~\eqref{eq-sitau} and the def\/inition of the root sequence, one obtains that
\begin{gather*}
\beta_{2k+1} =\eta(c_1) \eta(c_2)\cdots\eta(c_{2k-1}) \eta(c_{2k}) (\alpha_1)
\\
\phantom{\beta_{2k+1}}
 =\eta\big({-}a_{12}^X\big) \eta\big({-}a_{21}^{r_1(X)}\big)\cdots\eta\big({-}a_{12}^{(r_2r_1)^{k-1}(X)}\big)
 \eta\big({-}a_{21}^{r_1(r_2r_1)^{k-1}(X)}\big)(\alpha_1)
\\
\phantom{\beta_{2k+1}}
 =\eta\big({-}a_{12}^X\big)\tau \tau \eta\big({-}a_{21}^{r_1(X)}\big) \cdots\eta\big({-}a_{12}^{(r_2r_1)^{k-1}(X)}\big)\tau \tau
\eta\big({-}a_{21}^{r_1(r_2r_1)^{k-1}(X)}\big) (\alpha_1)
\\
\phantom{\beta_{2k+1}}
 =s_1^X s_2^{r_1(X)}\cdots s_1^{(r_2r_1)^{k-1}(X)}s_2^{r_1(r_2r_1)^{k-1}(X)}(\alpha_1)
 =\mathrm{id}_{X}(s_1s_2)^k \alpha_1,
\\
\tau \gamma_{2k+1} =\tau \eta\big({-}a_{21}^X\big) \eta\big({-}a_{12}^{r_2(X)}\big)\cdots \eta\big({-}a_{21}^{(r_1r_2)^{n-2}(X)}\big)
\eta\big({-}a_{12}^{r_2(r_1r_2)^{k-1}(X)}\big) (\alpha_1)
\\
\phantom{\tau \gamma_{2k+1}}
 =\big(\tau \eta\big({-}a_{21}^X\big) \eta\big({-}a_{12}^{r_2(X)}\big)\tau\big) (\tau \cdots \tau)\big(\tau \eta\big({-}a_{21}^{(r_1r_2)^{k-2}(X)}\big)
\eta\big({-}a_{12}^{r_2(r_1r_2)^{k-1}(X)}\big)\tau\big) \tau (\alpha_1)
\\
\phantom{\tau \gamma_{2k+1}}
 =\mathrm{id}_X(s_2s_1)^k \alpha_2.
\end{gather*}

The claims $\beta_{2k+2}=\mathrm{id}_X(s_1s_2)^k s_1 \alpha_2, \tau \gamma_{2k+2}=\mathrm{id}_X(s_2 s_1)^k s_2
\alpha_1$ hold by a~similar argument.

Thus $\rersys{X}=\{\pm \beta_k, \pm\tau\gamma_k|k\geq 1\}$ follows from the def\/inition of $\rersys{X}$.
\end{proof}

For a~f\/inite sequence $(v_1, \dots, v_n)$ of integers or vectors, where $n\geq 1$, let $(v_1, \dots,
v_n)^{\infty}=(u_k)_{k\geq1}$ be the sequence where $u_{mn+i}=v_i$ for all $1\leq i\leq n$, $m\geq0$.
\begin{Theorem}
\label{theo.cartan}
Let $\mathcal{C}= \mathcal{C}(I=\{1, 2\}, \mathcal{X}, r, (A^X)_{X\in \mathcal{X}})$ be a~connected semi-Cartan graph of rank two such that $|\mathcal{X}|$ is
finite.
Let $X\in \mathcal{X}$ and let~$n$ be the smallest positive integer with $(r_2r_1)^n(X)=X$.
Let $(c_k)_{k\geq 1}$ be the characteristic sequence of $\mathcal{C}$ with respect to~$X$ and~$1$, and let
$l=6n-\sum\limits_{i=1}^{2n}c_i$.
The following are equivalent:
\begin{itemize}\itemsep=0pt
\item[$(1)$] $\mathcal{C}$ is a~finite Cartan graph,
\item[$(2)$] $l>0$, $l |  12$, $(c_1, c_2, \dots,c_{12n/l})\in \mathcal{A}^{+}$, and $(c_k)_{k\geq1}=(c_1, c_2,
\dots,c_{12n/l})^{\infty}$.
\end{itemize}
In this case $12n/ l=|\rersys{X}_{+}|=t_{12}^{X}$.
\end{Theorem}

\begin{proof}
Let $(\beta_k)_{k\geq1}$ be the root sequence of $\mathcal{C}$ with respect to~$X$ and $1$ and $(\gamma_k)_{k\geq1}$ the root
sequence of $\mathcal{C}$ with respect to~$X$ and~$2$.

(1)$\Rightarrow$(2).
Let $q=t_{12}^X$.
Then $\rersys{Y}\subset \mathbb{N}_{0}^2\cup -\mathbb{N}_{0}^2$ for all $Y\in \mathcal{X}$ since $\mathcal{C}$ is a~Cartan graph.
By~\cite[Lemmas~3,~4]{HS10} and Lemma~\ref{le:root}, we have $\beta_k\in \mathbb{N}_0^2$ for all $1\leq k\leq q$ and
$\beta_{q}=\eta(c_1)\cdots\eta(c_{q-1})(\alpha_1)=\alpha_2$.
By the same reason we obtain that $\eta(c_2)\cdots\eta(c_{q})(\alpha_1)=\alpha_2$.
Then we have
\begin{gather*}
-\beta_{q+1}=-\eta(c_1)\cdots\eta(c_{q})(\alpha_1)=-\eta(c_1)(\alpha_2)=\alpha_1.
\end{gather*}
Thus $-\eta(c_1)\cdots\eta(c_{q})=\mathrm{id}$.
Indeed, if we set $w:=-\eta(c_1)\dots\eta(c_{q})$ and $w(\alpha_2):=a\alpha_1+b\alpha_2$, then $b=1$ since
$\det(w)=1$.
If~$q$ is odd then $-w\tau \in \Hom(Y, X)$ by equation~\eqref{eq-etadefi} and equation~\eqref{eq-sitau}, where
$Y=r_1(r_2r_1)^{(q-1)/2}(X)$.
In the same way, one gets $-w\in \Hom((r_2r_1)^{q/2}(X), X)$ if~$q$ is even.
Hence $w(\alpha_1), w(\alpha_2)\in \rersys{X}$.
Then $a\geq0$ since $w(\alpha_2)=\eta(c_1)\cdots\eta(c_{q-1})(\alpha_1)=\beta_q\in \rersys{X}\subset \mathbb{N}_{0}^2\cup
-\mathbb{N}_{0}^2$.
Moreover, $a\leq0$ since $w^{-1}(\alpha_2)=\alpha_2-a\alpha_1\in \mathbb{N}_{0}^2\cup -\mathbb{N}_{0}^2$.
Hence $a=0$ and $w(\alpha_2)=\alpha_2$.
Then $-\eta(c_1)\cdots\eta(c_{q})=\mathrm{id}$.
Hence $(c_1, \dots, c_q)\in \mathcal{A}^+$ by Theorem~\ref{theo:cqid}.
Therefore $\sum\limits_{i=1}^q c_i=3q-6$ by Lemma~\ref{le:sumA}.

Further, we apply the f\/irst part of the proof for $r_1(X)$ and the label $2$ instead of~$X$ and the label~$1$,
respectively.
Then $(c_2, \dots, c_{q+1})\in \mathcal{A}^+$ by Remark~\ref{re:char_se_a_j} and $\sum\limits_{i=2}^{q+1} c_i=3q-6$.
Hence \mbox{$c_{q+1}=c_1$}.
By induction, $(c_k)_{k\geq1}=(c_1, c_2, \dots,c_{q})^{\infty}$.
In particular, we obtain that $\sum\limits_{i=1}^{2qn} c_i=2n(3q-6)=q\big(\sum\limits_{i=1}^{2n} c_i\big)$.
Therefore $\sum\limits_{i=1}^{2n} c_i=6n-12n/q$.
Hence $q|12n$ and $l=6n-\sum\limits_{i=1}^{2n} c_i=12n/q>0$.
Further, $n|q$ since $\mathcal{C}$ is a~Cartan graph.
Thus $l|12$.

(2)$\Rightarrow$(1).
Set $q=12n/l$.
Then $\eta(c_1)\dots\eta(c_{q})=-\mathrm{id}$ and $\beta_k\in \mathbb{N}_0^2$ for $1\leq k\leq q$ by Theorem~\ref{theo:cqid}.
Then $(c_q, c_{q-1}, \dots,c_{1})\in \mathcal{A}^{+}$ by Corollary~\ref{coro:A+}(1) since $(c_1, c_2, \dots, c_{q})\in \mathcal{A}^{+}$.
Since $l|12$,~$q$ is a~multiple of~$n$.
Hence $(c_q, c_{q-1}, \dots,c_{1})^{\infty}$ is the characteristic sequence of $\mathcal{C}$ with respect to~$X$ and~$2$.
By Lemma~\ref{le:root}, we get that $\gamma_k \in \mathbb{N}_0^2$ for all $1\leq k\leq q$.
Therefore, since $(c_k)_{k\geq1}=(c_1, \dots, c_q)^{\infty}$, $\rersys{X}=\{\pm\beta_k, \pm\tau\gamma_k\,|\, 1\leq k\leq
q\}\subseteq \mathbb{N}_0^2 \cup -\mathbb{N}_0^2$ for all $X\in \mathcal{X}$ by Lemma~\ref{le:root}.
Hence $\mathcal{C}$ is f\/inite.

By the def\/inition of $t_{12}^X$ and~\cite[Lemma 4]{HS10}, we obtain that $t_{12}^X=q=|\rersys{X}_{+}|$.
Hence $n|t_{12}^X$ by assumption and $(r_2r_1)^{t_{12}^X}(X)=X$.
Therefore, $\mathcal{C}$ is a~Cartan graph.
\end{proof}

\section{Classif\/ication of rank two Nichols algebras of diagonal type\\
over f\/ields of positive characteristic}
\label{se:cla_Nichol}

{\samepage In this section, we classify all the two-dimensional braided vector spaces~$V$ of diagonal type over f\/ields of positive
characteristic such that the Nichols algebra of~$V$ has a~f\/inite root system.
The proof uses the characterization of the f\/inite Cartan graphs of rank two.

}

Let~$V$ be a~braided vector space of diagonal type with a~basis $\{x_1, x_2\}$ and the braiding $c (x_i\otimes
x_j)=q_{ij}x_j\otimes x_i$, where $q_{ij}\in \Bbbk$, $i,j\in \{1, 2\}$.
We choose an abelian group~$G$ and the set $\{g_i\,|\, g_i\in G,\, i\in \{1, 2\}\}$ such that the assignments $\delta
(x_i)=g_i x_i$, $g_i.
x_j=q_{ij}x_j$ for $i,j\in \{1, 2\}$ def\/ine a~Yetter--Drinfel'd module structure on~$V$ over~$\Bbbk G$.
Let $\mathcal{B}(V)$ denote the Nichols algebra of~$V$.

The following theorem determines whether Weyl groupoid $\mathcal{W}(\Bbbk x_1,\Bbbk x_2)$ of $(\Bbbk x_1,\Bbbk x_2)$ is f\/inite
in terms of the Dynkin diagram of~$V$.
\begin{Theorem}
\label{Theo:clasi}
Let~$V$ be a~two-dimensional braided vector space of diagonal type with the brai\-ding
\begin{gather*}
c (x_i\otimes x_j)=q_{ij}x_j\otimes x_i,
\end{gather*}
where $i, j\in \{1, 2\}$ and $\{x_1, x_2\}$ is a~basis of~$V$.
Let $M:=(\Bbbk x_1,\Bbbk x_2)$.
Assume that the characteristic~$p$ of~$\Bbbk$ is positive.
Then the following are equivalent:
\begin{enumerate}\itemsep=0pt
\item[$(1)$] $\mathcal{B}(V)$ is decomposable and ${\boldsymbol{\Delta}}^{[M]}$ is finite,
\item[$(2)$] the Dynkin diagram $\mathcal{D}$ of~$V$ appears in Tables~{\rm \ref{tab.1}},~{\rm \ref{tab.2}},~{\rm \ref{tab.3}},~{\rm \ref{tab.4}}
and~{\rm \ref{tab.5}}, if $p=2$, $p=3$, $p=5$, $p=7$ and $p>7$, respectively,
\item[$(3)$]~$M$ admits all refections and $\mathcal{W}(M)$ is finite.
\end{enumerate}
In this case, the row of Table~{\rm \ref{tab.6}} containing~$\mathcal{D}$ consists precisely of the Dynkin diagrams of the points of~$\mathcal{C}_s(M)$.
Further, the corresponding row of Table~{\rm \ref{tab.7}} contains the exchange graph of~$\mathcal{C}_s(M)$.
\end{Theorem}

{\setlength{\unitlength}{1mm} \settowidth{\mpb}{$q_0\in \Bbbk^\ast \setminus \{-1,1\}$,}
\begin{table}[t]\centering
\caption{Dynkin diagrams in characteristic $p=2$.}\label{tab.1}
\vspace{1mm}
\begin{tabular}{r|l|l}
\hline
& \text{Dynkin diagrams} & \text{f\/ixed parameters}
\\
\hline
\hline
$1$ & \rule[-3\unitlength]{0pt}{8\unitlength}
\begin{picture}
(14,5)(0,3) \put(1,2){\circle{2}} \put(13,2){\circle{2}} \put(1,5){\makebox[0pt]{\scriptsize~$q$}}
\put(13,5){\makebox[0pt]{\scriptsize~$r$}}
\end{picture}
& $q,r\in \Bbbk^\ast $
\\
\hline
$2$ & \Dchaintwo{}{$q$}{$q^{-1}$}{$q$} & $q\in \Bbbk^\ast \setminus \{1\}$
\\
\hline
$3$ & \Dchaintwo{}{$q$}{$q^{-1}$}{$1$}\ \Dchaintwo{}{$1$}{$q$}{$1$} & $q\in \Bbbk^\ast \setminus \{1\}$
\\
\hline
$4$ & \Dchaintwo{}{$q$}{$q^{-2}$}{$q^2$} & $q\in \Bbbk^\ast \setminus \{1\}$
\\
\hline
\rule[-3\unitlength]{0pt}{10\unitlength} $5$ & \Dchaintwo{}{$q$}{$q^{-2}$}{$1$}\ \Dchaintwo{}{$q^{-1}$}{$q^2$}{$1$} &
$q\in \Bbbk^\ast \setminus \{1\}$
\\
\hline
\rule[-3\unitlength]{0pt}{10\unitlength} \settowidth{\mpb}{$q_0\in \Bbbk^\ast \setminus \{-1,1\}$,} 6 &
\Dchaintwo{}{$\zeta $}{$q^{-1}$}{$q$}\ \Dchaintwo{}{$\zeta $}{$\zeta^{-1}q$}{\ \ $\zeta q^{-1}$} & $\zeta \in G'_3$,
$q\in \Bbbk^\ast \setminus \{1,\zeta,\zeta^2\}$
\\
\hline
$7$ & \Dchaintwo{}{$\zeta $}{$\zeta $}{$1$}\ \Dchaintwo{}{$\zeta^{-1}$}{$\zeta^{-1}$}{$1$} & $\zeta \in G'_3$
\\
\hline
$10$ & \Dchaintwo{}{$\zeta^2$}{$\zeta $}{$1$}\ \Dchaintwo{}{$\zeta^3$}{$\zeta^{-1}$}{$1$}\ \Dchaintwo{}{$\zeta^3
$}{$\zeta^{-2}$}{$\zeta $} & $\zeta \in G'_9$
\\
\hline
\rule[-3\unitlength]{0pt}{10\unitlength} $11$ & \Dchaintwo{}{$q$}{$q^{-3}$}{$q^3$} & $q\in \Bbbk^\ast \setminus \{1\}$,
$q\notin G'_3$
\\
\hline
$14$ & \Dchaintwo{}{$\zeta $}{$\zeta^2$}{$1$}\ \Dchaintwo{}{$\zeta^{-2}$}{$\zeta^{-2}$}{$1$} & $\zeta \in G'_5$
\\
\hline
$16$ & \Dchaintwo{}{$\zeta^5 $}{$\zeta^{-3}$}{$\zeta $}\ \Dchaintwo{}{$\zeta^5$}{$\zeta^{-2}$}{$1$}\
\Dchaintwo{}{$\zeta^3$}{$\zeta^2$}{$1$}\ \Dchaintwo{}{$\zeta^3$}{$\zeta^4$}{$\zeta^{-4}$} & $\zeta \in G'_{15}$
\\
\hline
$17$ & \Dchaintwo{}{$\zeta $}{$\zeta^{-3}$}{$1$}\ \Dchaintwo{}{$\zeta^{-2}$}{$\zeta^3$}{$1$} & $\zeta \in G'_7$
\\
\hline
\end{tabular}
\end{table}}

{\setlength{\unitlength}{1mm}
\begin{table}[t]\centering
\caption{Dynkin diagrams in characteristic $p=3$.}\label{tab.2}
\vspace{1mm}

\begin{tabular}{r|l|l}
\hline
& \text{Dynkin diagrams} & \text{f\/ixed parameters}
\\
\hline
\hline
1 & \rule[-3\unitlength]{0pt}{8\unitlength}
\begin{picture}
(14,5)(0,3) \put(1,2){\circle{2}} \put(13,2){\circle{2}} \put(1,5){\makebox[0pt]{\scriptsize~$q$}}
\put(13,5){\makebox[0pt]{\scriptsize~$r$}}
\end{picture}
& $q,r\in \Bbbk^\ast $
\\
\hline
2 & \Dchaintwo{}{$q$}{$q^{-1}$}{$q$} & $q\in \Bbbk^\ast \setminus \{1\}$
\\
\hline
3 & \Dchaintwo{}{$q$}{$q^{-1}$}{$-1$}\ \Dchaintwo{}{$-1$}{$q$}{$-1$} & $q\in \Bbbk^\ast \setminus \{-1,1\}$
\\
\hline
4 & \Dchaintwo{}{$q$}{$q^{-2}$}{$q^2$} & $q\in \Bbbk^\ast \setminus \{-1,1\}$
\\
\hline
\rule[-3\unitlength]{0pt}{10\unitlength} 5 & \Dchaintwo{}{$q$}{$q^{-2}$}{$-1$}\ \Dchaintwo{}{$-q^{-1}$}{$q^2$}{$-1$} &
$q\in \Bbbk^\ast \setminus \{-1,1\}$, $q\notin G'_4$
\\
\hline
\rule[-3\unitlength]{0pt}{10\unitlength} $6' $& \Dchaintwo{}{$1$}{$q$}{$q^{-1}$}\ \Dchaintwo{}{$1$}{$q^{-1}$}{$q$} &
$q\in \Bbbk^\ast \setminus \{1, -1\}$
\\
\hline
$6''' $& \Dchaintwo{}{$1$}{$-1$}{$-1$} &
\\
\hline
$9'$& \Dchaintwo{}{$\zeta $}{$\zeta$}{$-1$}\ \Dchaintwo{}{$1$}{$-\zeta$}{$-1$}\ \Dchaintwo{}{$1$}{$\zeta $}{$1$} &
$\zeta \in G'_{4}$
\\
\hline
11& \Dchaintwo{}{$q$}{$q^{-3}$}{$q^3$} & $q\in \Bbbk^\ast \setminus \{-1,1\}$
\\
\hline
12 & \Dchaintwo{}{$\zeta $}{$-\zeta $}{$-1$}\ \Dchaintwo{}{$\zeta^2$}{$-\zeta^{-1}$}{$-1$}\
\Dchaintwo{}{$\zeta^2$}{$\zeta $}{$\zeta^{-1}$} & $\zeta \in G'_8$
\\
\hline
$13'$& \Dchaintwo{}{$-\zeta $}{$\zeta^{-1}$}{$-1$}\ \Dchaintwo{}{$1$}{$\zeta$}{$-1$}\
\Dchaintwo{}{$1$}{$\zeta^{-1}$}{$-\zeta^{2}$}\ \ \Dchaintwo{}{\ \ $-\zeta^{-1}$}{$-\zeta $}{$-\zeta^2$}& $\zeta \in G'_{8}$
\\
\hline
14 & \Dchaintwo{}{$\zeta $}{$\zeta^2$}{$-1$}\ \Dchaintwo{}{$-\zeta^{-2}$}{$\zeta^{-2}$}{$-1$} & $\zeta \in G'_5$
\\
\hline
15 & \Dchaintwo{}{$\zeta $}{$\zeta^{-3}$}{$-1$}\ \Dchaintwo{}{$-\zeta^{-2}$}{$\zeta^3$}{$-1$}\
\Dchaintwo{}{$-\zeta^{-2}$}{$-\zeta^3$}{$-1$}\ \Dchaintwo{}{$-\zeta $}{$-\zeta^{-3}$}{$-1$} & $\zeta \in G'_{20}$
\\
\hline
$16'$ & \Dchaintwo{}{$1$}{$-\zeta^{-1}$}{$-\zeta^2 $}\ \Dchaintwo{}{$1$}{$-\zeta$}{$-1$}\
\Dchaintwo{}{$\zeta$}{$-\zeta^{-1}$}{$-1$}\ \Dchaintwo{}{$\zeta$}{$-\zeta^3$}{$-\zeta^{-3}$} & $\zeta \in G'_{5}$
\\
\hline
17 & \Dchaintwo{}{$-\zeta $}{$-\zeta^{-3}$}{$-1$}\ \Dchaintwo{}{$-\zeta^{-2}$}{$-\zeta^3$}{$-1$} & $\zeta \in G'_7$
\\
\hline
\end{tabular}
\end{table}}

{\setlength{\unitlength}{1mm}
\begin{table}[t]\centering
\caption{Dynkin diagrams in characteristic $p=5$.}\label{tab.3}
\vspace{1mm}
\begin{tabular}{r|l|l}
\hline
& \text{Dynkin diagrams} & \text{f\/ixed parameters}
\\
\hline
\hline
1 & \rule[-3\unitlength]{0pt}{8\unitlength}
\begin{picture}
(14,5)(0,3) \put(1,2){\circle{2}} \put(13,2){\circle{2}} \put(1,5){\makebox[0pt]{\scriptsize~$q$}}
\put(13,5){\makebox[0pt]{\scriptsize~$r$}}
\end{picture}
& $q,r\in \Bbbk^\ast $
\\
\hline
2 & \Dchaintwo{}{$q$}{$q^{-1}$}{$q$} & $q\in \Bbbk^\ast \setminus \{1\}$
\\
\hline
3 & \Dchaintwo{}{$q$}{$q^{-1}$}{$-1$}\ \Dchaintwo{}{$-1$}{$q$}{$-1$} & $q\in \Bbbk^\ast \setminus \{-1,1\}$
\\
\hline
4 & \Dchaintwo{}{$q$}{$q^{-2}$}{$q^2$} & $q\in \Bbbk^\ast \setminus \{-1,1\}$
\\
\hline
\rule[-3\unitlength]{0pt}{10\unitlength} 5 & \Dchaintwo{}{$q$}{$q^{-2}$}{$-1$}\ \Dchaintwo{}{$-q^{-1}$}{$q^2$}{$-1$} &
$q\in \Bbbk^\ast \setminus \{-1,1\}$, $q\notin G'_4$
\\
\hline
\rule[-3\unitlength]{0pt}{10\unitlength} \settowidth{\mpb}{$q_0\in \Bbbk^\ast \setminus \{-1,1\}$,} 6 &
\Dchaintwo{}{$\zeta $}{$q^{-1}$}{$q$}\ \Dchaintwo{}{$\zeta $}{$\zeta^{-1}q$}{\ \ $\zeta q^{-1}$} & $\begin{array}{@{}c@{}}\zeta
\in G'_3, \ q \zeta\not=-1 \\
q\in \Bbbk^\ast \setminus \{1,\zeta,\zeta^2\}\end{array}$
\\
\hline
$ 6''$ & \Dchaintwo{}{$\zeta$}{$-\zeta$}{$-{\zeta}^{-1}$} & $\zeta \in G'_3$
\\
\hline
7 & \Dchaintwo{}{$\zeta $}{$-\zeta $}{$-1$}\ \Dchaintwo{}{$\zeta^{-1}$}{$-\zeta^{-1}$}{$-1$} & $\zeta \in G'_3$
\\
\hline
8 &\Dchaintwo{}{$-\zeta^{2}$}{$\zeta$}{$-1$}\ \Dchaintwo{}{$-\zeta^2$}{$\zeta^3$}{$-\zeta^{-2}$}\
\Dchaintwo{}{$-1$}{$-\zeta^{-1}$}{$\ \ -\zeta^{-2}$}\  \ \ \Dchaintwo{}{$-1$}{$-\zeta$}{$\zeta^3$}\
\Dchaintwo{}{$-1$}{$\zeta^{-1}$}{$\zeta^3$} & $\zeta \in G'_{12}$
\\
\hline
9 & \Dchaintwo{}{$-\zeta^{-1}$}{$-\zeta^3$}{$-1$}\ \Dchaintwo{}{$-\zeta^2$}{$\zeta^3$}{$-1$}\
\Dchaintwo{}{$-\zeta^2$}{$\zeta $}{$-\zeta^2$} & $\zeta \in G'_{12}$
\\
\hline
10 & \Dchaintwo{}{$-\zeta^2$}{$\zeta $}{$-1$}\ \Dchaintwo{}{$\zeta^3$}{$\zeta^{-1}$}{$-1$}\
\Dchaintwo{}{$\zeta^3$}{$\zeta^{-2}$}{$-\zeta $} & $\zeta \in G'_9$
\\
\hline
\rule[-3\unitlength]{0pt}{10\unitlength} 11 & \Dchaintwo{}{$q$}{$q^{-3}$}{$q^3$} & $q\in \Bbbk^\ast \setminus \{-1,1\}$,
$q\notin G'_3$
\\
\hline
12 & \Dchaintwo{}{$\zeta $}{$-\zeta $}{$-1$}\ \Dchaintwo{}{$\zeta^2$}{$-\zeta^{-1}$}{$-1$}\
\Dchaintwo{}{$\zeta^2$}{$\zeta $}{$\zeta^{-1}$} & $\zeta \in G'_8$
\\
\hline
13 & \Dchaintwo{}{$\zeta $}{$\zeta^{-5}$}{$-1$}\ \Dchaintwo{}{\ \ $-\zeta^{-4}$}{\ \ $\zeta^5$}{$-1$}\
\Dchaintwo{}{$-\zeta^{-4}$}{$-\zeta^{-1}$}{$\zeta^6$}\ \Dchaintwo{}{$\zeta^{-1}$}{$\zeta $}{$\zeta^6$} & $\zeta \in
G'_{24}$
\\
\hline
$15'$ & \Dchaintwo{}{$\zeta $}{$\zeta $}{$-1$}\ \Dchaintwo{}{$1$}{$-\zeta $}{$-1$}\ \Dchaintwo{}{$1$}{$\zeta$}{$-1$}\
\Dchaintwo{}{$-\zeta $}{$-\zeta $}{$-1$}\ & $\zeta \in G'_{4}$
\\
\hline
$16''$ & \Dchaintwo{}{$\zeta^2$}{$-1$}{$-\zeta $}\ \Dchaintwo{}{$\zeta^2$}{$-\zeta $}{$-1$}\
\Dchaintwo{}{$1$}{$-\zeta^{-1}$}{$-1$}\ \Dchaintwo{}{$1$}{$-\zeta$}{$-\zeta^{-1}$} & $\zeta \in G'_{3}$
\\
\hline
17 & \Dchaintwo{}{$-\zeta $}{$-\zeta^{-3}$}{$-1$}\ \Dchaintwo{}{\ \ $-\zeta^{-2}$}{\ \ $-\zeta^3$}{$-1$} & $\zeta \in
G'_7$
\\
\hline
\end{tabular}
\end{table}}

{\setlength{\unitlength}{1mm}
\begin{table}[t]\centering
\caption{Dynkin diagrams in characteristic $p=7$.}\label{tab.4}
\vspace{1mm}
\begin{tabular}{r|l|l}
\hline
& \text{Dynkin diagrams} & \text{f\/ixed parameters}
\\
\hline
\hline
1 & \rule[-3\unitlength]{0pt}{8\unitlength}
\begin{picture}
(14,5)(0,3) \put(1,2){\circle{2}} \put(13,2){\circle{2}} \put(1,5){\makebox[0pt]{\scriptsize~$q$}}
\put(13,5){\makebox[0pt]{\scriptsize~$r$}}
\end{picture}
& $q,r\in \Bbbk^\ast $
\\
\hline
2 & \Dchaintwo{}{$q$}{$q^{-1}$}{$q$} & $q\in \Bbbk^\ast \setminus \{1\}$
\\
\hline
3 & \Dchaintwo{}{$q$}{$q^{-1}$}{$-1$}\ \Dchaintwo{}{$-1$}{$q$}{$-1$} & $q\in \Bbbk^\ast \setminus \{-1,1\}$
\\
\hline
4 & \Dchaintwo{}{$q$}{$q^{-2}$}{$q^2$} & $q\in \Bbbk^\ast \setminus \{-1,1\}$
\\
\hline
\rule[-3\unitlength]{0pt}{10\unitlength} 5 & \Dchaintwo{}{$q$}{$q^{-2}$}{$-1$}\ \Dchaintwo{}{$-q^{-1}$}{$q^2$}{$-1$} &
$q\in \Bbbk^\ast \setminus \{-1,1\}$, $q\notin G'_4$
\\
\hline
\rule[-3\unitlength]{0pt}{10\unitlength} \settowidth{\mpb}{$q_0\in \Bbbk^\ast \setminus \{-1,1\}$,} 6 &
\Dchaintwo{}{$\zeta $}{$q^{-1}$}{$q$}\ \Dchaintwo{}{$\zeta $}{$\zeta^{-1}q$}{\ \ $\zeta q^{-1}$} &
$\begin{array}{@{}l@{}}
\zeta \in G'_3, \ q \zeta\not=-1 \\
q\in \Bbbk^\ast \setminus \{1,\zeta,\zeta^2\}
\end{array}$
\\
\hline
$6''$ & \Dchaintwo{}{$\zeta$}{$-\zeta$}{$-{\zeta}^{-1}$} & $\zeta \in G'_3$
\\
\hline
7 & \Dchaintwo{}{$\zeta $}{$-\zeta $}{$-1$}\ \Dchaintwo{}{$\zeta^{-1}$}{$-\zeta^{-1}$}{$-1$} & $\zeta \in G'_3$
\\
\hline
8 &\Dchaintwo{}{$-\zeta^{2}$}{$\zeta$}{$-1$}\ \Dchaintwo{}{$-\zeta^2$}{$\zeta^3$}{$-\zeta^{-2}$}\
\Dchaintwo{}{$-1$}{$-\zeta^{-1}$}{$-\zeta^{-2}$}\  \ \ \ \Dchaintwo{}{$-1$}{$-\zeta$}{$\zeta^3$}\ \Dchaintwo{}{$-1$}{$\zeta^{-1}$}{$\zeta^3$} & $\zeta \in G'_{12}$
\\
\hline
9 &\Dchaintwo{}{$-\zeta^{-1}$}{$-\zeta^3$}{$-1$}\ \Dchaintwo{}{$-\zeta^2$}{$\zeta^3$}{$-1$}\
\Dchaintwo{}{$-\zeta^2$}{$\zeta $}{$-\zeta^2$} & $\zeta \in G'_{12}$
\\
\hline
10 & \Dchaintwo{}{$-\zeta^2$}{$\zeta $}{$-1$}\ \Dchaintwo{}{$\zeta^3$}{$\zeta^{-1}$}{$-1$}\
\Dchaintwo{}{$\zeta^3$}{$\zeta^{-2}$}{$-\zeta $} & $\zeta \in G'_9$
\\
\hline
\rule[-3\unitlength]{0pt}{10\unitlength} 11 & \Dchaintwo{}{$q$}{$q^{-3}$}{$q^3$} & $q\in \Bbbk^\ast \setminus \{-1,1\}$,
$q\notin G'_3$
\\
\hline
12 & \Dchaintwo{}{$\zeta $}{$-\zeta $}{$-1$}\ \Dchaintwo{}{$\zeta^2$}{$-\zeta^{-1}$}{$-1$}\
\Dchaintwo{}{$\zeta^2$}{$\zeta $}{$\zeta^{-1}$} & $\zeta \in G'_8$
\\
\hline
13 & \Dchaintwo{}{$\zeta $}{$\zeta^{-5}$}{$-1$}\ \Dchaintwo{}{$-\zeta^{-4}$}{$\zeta^5$}{$-1$}\
\Dchaintwo{}{$-\zeta^{-4}$}{$-\zeta^{-1}$}{$\zeta^6$}\ \Dchaintwo{}{$\zeta^{-1}$}{$\zeta $}{$\zeta^6$} & $\zeta \in G'_{24}$
\\
\hline
14 & \Dchaintwo{}{$\zeta $}{$\zeta^2$}{$-1$}\ \Dchaintwo{}{$-\zeta^{-2}$}{$\zeta^{-2}$}{$-1$} & $\zeta \in G'_5$
\\
\hline
15 & \Dchaintwo{}{$\zeta $}{$\zeta^{-3}$}{$-1$}\ \Dchaintwo{}{$-\zeta^{-2}$}{$\zeta^3$}{$-1$}\
\Dchaintwo{}{$-\zeta^{-2}$}{$-\zeta^3$}{$-1$}\ \Dchaintwo{}{$-\zeta $}{$-\zeta^{-3}$}{$-1$} & $\zeta \in G'_{20}$
\\
\hline
16 & \Dchaintwo{}{$\zeta^5 $}{$-\zeta^{-3}$}{$-\zeta$}\ \Dchaintwo{}{$\zeta^5$}{$-\zeta^{-2}$}{$-1$}\
\Dchaintwo{}{$\zeta^3$}{$-\zeta^2$}{$-1$}\ \Dchaintwo{}{$\zeta^3$}{$-\zeta^4$}{$-\zeta^{-4}$} & $\zeta \in G'_{15}$
\\
\hline
18 & \Dchaintwo{}{$\zeta^{-1} $}{$-1$}{$-\zeta$}\ \Dchaintwo{}{$\zeta^{-1}$}{$-\zeta$}{$-1$}\ \Dchaintwo{}{$1$}{$-\zeta^{-1}$}{$-1$} & $\zeta \in G'_3$
\\
& \Dchaintwo{}{$1$}{$-\zeta$}{$-1$}\ \Dchaintwo{}{$\zeta$}{$-\zeta^{-1}$}{$-1$}\
\Dchaintwo{}{$\zeta$}{$-1$}{$-\zeta^{-1} $} &
\\
\hline
\end{tabular}
\end{table}}

{\setlength{\unitlength}{1mm} \settowidth{\mpb}{$q_0\in \Bbbk^\ast \setminus \{-1,1\}$,}

\begin{table}[t]\centering
\caption{Dynkin diagrams in characteristic $p>7$.} \label{tab.5}
\vspace{1mm}
\begin{tabular}{r|l|l}
\hline
& \text{Dynkin diagrams} & \text{f\/ixed parameters}
\\
\hline
\hline
1 & \rule[-3\unitlength]{0pt}{8\unitlength}
\begin{picture}
(14,5)(0,3) \put(1,2){\circle{2}} \put(13,2){\circle{2}} \put(1,5){\makebox[0pt]{\scriptsize~$q$}}
\put(13,5){\makebox[0pt]{\scriptsize~$r$}}
\end{picture}
& $q,r\in \Bbbk^\ast $
\\
\hline
2 & \Dchaintwo{}{$q$}{$q^{-1}$}{$q$} & $q\in \Bbbk^\ast \setminus \{1\}$
\\
\hline
3 & \Dchaintwo{}{$q$}{$q^{-1}$}{$-1$}\ \Dchaintwo{}{$-1$}{$q$}{$-1$} & $q\in \Bbbk^\ast \setminus \{-1,1\}$
\\
\hline
4 & \Dchaintwo{}{$q$}{$q^{-2}$}{$q^2$} & $q\in \Bbbk^\ast \setminus \{-1,1\}$
\\
\hline
\rule[-3\unitlength]{0pt}{8\unitlength} 5 & \Dchaintwo{}{$q$}{$q^{-2}$}{$-1$}\ \Dchaintwo{}{
$-q^{-1}$}{$q^2$}{$-1$} & \parbox{\mpb}{$q\in \Bbbk^\ast \setminus \{-1,1\}$
\\
$q\notin G'_4$}
\\
\hline
\rule[-3\unitlength]{0pt}{10\unitlength} \settowidth{\mpb}{$q_0\in \Bbbk^\ast \setminus \{-1,1\}$,} 6 &
\Dchaintwo{}{$\zeta $}{$q^{-1}$}{$q$}\ \Dchaintwo{}{$\zeta $}{$\zeta^{-1}q$}{\ \ $\zeta q^{-1}$} & \parbox{\mpb}{$\zeta
\in G'_3$, $q \zeta\not=-1$
\\
$q\in \Bbbk^\ast \setminus \{1,\zeta,\zeta^2\}$}
\\
\hline
$6''$ & \Dchaintwo{}{$\zeta$}{$\!\!-\zeta$}{$\ -{\zeta}^{-1}$} & $\zeta \in G'_3$
\\
\hline
7 & \Dchaintwo{}{$\zeta $}{$-\zeta $}{$-1$}\ \Dchaintwo{}{$\zeta^{-1}$}{$-\zeta^{-1}$}{$-1$} & $\zeta \in G'_3$
\\
\hline
8 &\Dchaintwo{}{$-\zeta^{2}$}{$\zeta$}{$-1$}\ \Dchaintwo{}{$-\zeta^2$}{$\zeta^3$}{$-\zeta^{-2}$}\
\Dchaintwo{}{$-1$}{$-\zeta^{-1}$}{\ \ $-\zeta^{-2}$}\ \ \ \Dchaintwo{}{$-1$}{$-\zeta$}{$\zeta^3$}\ \Dchaintwo{}{$-1$}{$\zeta^{-1}$}{$\zeta^3$} & $\zeta \in G'_{12}$
\\
\hline
9 & \Dchaintwo{}{$-\zeta^{-1}$}{$-\zeta^3$}{$-1$}\ \Dchaintwo{}{$-\zeta^2$}{$\zeta^3$}{$-1$}\
\Dchaintwo{}{$-\zeta^2$}{$\zeta $}{$-\zeta^2$} & $\zeta \in G'_{12}$
\\
\hline
10 & \Dchaintwo{}{$-\zeta^2$}{$\zeta $}{$-1$}\ \Dchaintwo{}{$\zeta^3$}{$\zeta^{-1}$}{$-1$}\ \Dchaintwo{}{$\zeta^3
$}{$\zeta^{-2}$}{$-\zeta$} & $\zeta \in G'_9$
\\
\hline
\rule[-3\unitlength]{0pt}{10\unitlength} 11 & \Dchaintwo{}{$q$}{$q^{-3}$}{$q^3$} & \parbox{\mpb}{$q\in \Bbbk^\ast
\setminus \{-1,1\}$
\\
$q\notin G'_3$}
\\
\hline
12 & \Dchaintwo{}{$\zeta $}{$-\zeta $}{$-1$}\ \Dchaintwo{}{$\zeta^2$}{$-\zeta^{-1}$}{$-1$}\
\Dchaintwo{}{$\zeta^2$}{$\zeta $}{$\zeta^{-1}$} & $\zeta \in G'_8$
\\
\hline
13 &\Dchaintwo{}{$\zeta $}{$\zeta^{-5}$}{$-1$}\ \Dchaintwo{}{$-\zeta^{-4}$}{$\zeta^5$}{$-1$}\
\Dchaintwo{}{$-\zeta^{-4}$}{$-\zeta^{-1}$}{$\zeta^6$}\ \Dchaintwo{}{$\zeta^{-1}$}{$\zeta $}{$\zeta^6$} & $\zeta \in G'_{24}$
\\
\hline
14 & \Dchaintwo{}{$\zeta $}{$\zeta^2$}{$-1$}\ \Dchaintwo{}{$-\zeta^{-2}$}{$\zeta^{-2}$}{$-1$} & $\zeta \in G'_5$
\\
\hline
15 & \Dchaintwo{}{$\zeta $}{$\zeta^{-3}$}{$-1$}\ \Dchaintwo{}{$-\zeta^{-2}$}{$\zeta^3$}{$-1$}\
\Dchaintwo{}{$-\zeta^{-2}$}{$-\zeta^3$}{$-1$}\ \Dchaintwo{}{$-\zeta $}{$-\zeta^{-3}$}{$-1$} & $\zeta \in G'_{20}$
\\
\hline
16 & \Dchaintwo{}{$\zeta^5 $}{$-\zeta^{-3}$}{$-\zeta$}\ \Dchaintwo{}{$\zeta^5$}{$-\zeta^{-2}$}{$-1$}\
\Dchaintwo{}{$\zeta^3$}{$-\zeta^2$}{$-1$}\ \Dchaintwo{}{$\zeta^3$}{$-\zeta^4$}{$-\zeta^{-4}$} & $\zeta \in G'_{15}$
\\
\hline
17 & \Dchaintwo{}{$-\zeta $}{$-\zeta^{-3}$}{$-1$}\ \Dchaintwo{}{$-\zeta^{-2}$}{$-\zeta^3$}{$-1$} & $\zeta \in G'_7$
\\
\hline
\end{tabular}
\end{table}}

{\setlength{\unitlength}{1mm} \settowidth{\mpb}{$q_0\in \Bbbk^\ast \setminus \{-1,1\}$,}
\begin{table}[t!] \centering
\caption{The Dynkin diagrams in Theorem~\ref{Theo:clasi}.} \label{tab.6}
\vspace{1mm}
\begin{tabular}{@{\,}r@{\,}|@{\,\,}l|@{\,}l@{\,}|@{\,}l@{\,}}
\hline
& \text{Dynkin diagrams} & \text{f\/ixed parameters} & \text{char} $\Bbbk$
\\
\hline
\hline
1 & \rule[-3\unitlength]{0pt}{8\unitlength}
\begin{picture}
(14,5)(0,3) \put(1,2){\circle{2}} \put(13,2){\circle{2}} \put(1,5){\makebox[0pt]{\scriptsize~$q$}}
\put(13,5){\makebox[0pt]{\scriptsize~$r$}}
\end{picture}
& $q,r\in \Bbbk^\ast $ &
\\
\hline
2 & \Dchaintwo{}{$q$}{$q^{-1}$}{$q$} & $q\in \Bbbk^\ast \setminus \{1\}$ &
\\
\hline
3 & \Dchaintwo{}{$q$}{$q^{-1}$}{$-1$}\ \Dchaintwo{}{$-1$}{$q$}{$-1$} & $q\in \Bbbk^\ast \setminus \{-1,1\}$ &
\\
\hline
4 & \Dchaintwo{}{$q$}{$q^{-2}$}{$q^2$} & $q\in \Bbbk^\ast \setminus \{-1,1\}$ &
\\
\hline
\rule[-3\unitlength]{0pt}{10\unitlength} 5 & \Dchaintwo{}{$q$}{$q^{-2}$}{$-1$}\ \Dchaintwo{}{$-q^{-1}$}{$q^2$}{$-1$} &
\parbox{\mpb}{$q\in \Bbbk^\ast \setminus \{-1,1\}$
\\
$q\notin G'_4$} &
\\
\hline
\rule[-3\unitlength]{0pt}{10\unitlength} 6 & \Dchaintwo{}{$\zeta $}{$q^{-1}$}{$q$}\ \Dchaintwo{}{$\zeta
$}{$\zeta^{-1}q$}{\ \ \ \ $\zeta q^{-1}$} & \parbox{\mpb}{$\zeta \in G'_3$, $q \zeta\not=-1$
\\
$q\in \Bbbk^\ast \setminus \{1,\zeta,\zeta^2\}$} & $p\not=3$
\\
\hline
$6'$ & \Dchaintwo{}{$1$}{$q$}{$q^{-1}$}\ \Dchaintwo{}{$1$}{$q^{-1}$}{$q$} & $q\in \Bbbk^\ast \setminus \{1, -1\}$ &
$p=3$
\\
\hline
$6''$ & \Dchaintwo{}{$\zeta$}{$-\zeta$}{$-{\zeta}^{-1}$} & $\zeta \in G'_3$ & $p\not=2, 3$
\\
\hline
$6'''$ & \Dchaintwo{}{$1$}{$-1$}{$-1$} & & $p=3$
\\
\hline
7 & \Dchaintwo{}{$\zeta $}{$-\zeta $}{$-1$}\ \Dchaintwo{}{$\zeta^{-1}$}{$-\zeta^{-1}$}{$-1$} & $\zeta \in G'_3$ & $p\not=3$
\\
\hline
8 &\Dchaintwo{}{$-\zeta^{2}$}{$\zeta$}{$-1$} \ \Dchaintwo{}{$-\zeta^2$}{$\zeta^3$}{$-\zeta^{-2}$} \
\Dchaintwo{}{$-1$}{$-\zeta^{-1}$}{ \ $-\zeta^{-2}$} \ \Dchaintwo{}{$-1$}{$-\zeta$}{$\zeta^3$}   \Dchaintwo{}{$-1$}{$\zeta^{-1}$}{$\zeta^3$} & $\zeta \in G'_{12}$ & $p\not=2, 3$
\\
\hline
9 & \Dchaintwo{}{$-\zeta^{-1}$}{$-\zeta^3$}{$-1$}\ \Dchaintwo{}{$-\zeta^2$}{$\zeta^3$}{$-1$}\
\Dchaintwo{}{$-\zeta^2$}{$\zeta $}{$-\zeta^2$} & $\zeta \in G'_{12}$ & $p\not=2,3 $
\\
\hline
$9'$& \Dchaintwo{}{$\zeta $}{$\zeta$}{$-1$}\ \Dchaintwo{}{$1$}{$-\zeta$}{$-1$}\ \Dchaintwo{}{$1$}{$\zeta $}{$1$} &
$\zeta \in G'_{4}$ & $p=3 $
\\
\hline
\rule[-3\unitlength]{0pt}{10\unitlength} 10 & \Dchaintwo{}{$-\zeta^2$}{$\zeta $}{$-1$}\
\Dchaintwo{}{$\zeta^3$}{$\zeta^{-1}$}{$-1$}\ \Dchaintwo{}{$\zeta^3 $}{$\zeta^{-2}$}{$-\zeta$} & $\zeta \in G'_9$ &
$p\not=3$
\\
\hline
11 & \Dchaintwo{}{$q$}{$q^{-3}$}{$q^3$} & \parbox{\mpb}{$q\in \Bbbk^\ast \setminus \{-1,1\}$
\\
$q\notin G'_3$} &
\\
\hline
12 & \Dchaintwo{}{$\zeta $}{$-\zeta $}{$-1$}\ \Dchaintwo{}{$\zeta^2$}{$-\zeta^{-1}$}{$-1$}\
\Dchaintwo{}{$\zeta^2$}{$\zeta $}{$\zeta^{-1}$} & $\zeta \in G'_8$ & $p\not=2$
\\
\hline
13 &\Dchaintwo{}{$\zeta $}{$\zeta^{-5}$}{$-1$}\ \Dchaintwo{}{$-\zeta^{-4}$}{$\zeta^5$}{$-1$}\
\Dchaintwo{}{$-\zeta^{-4}$}{$-\zeta^{-1}$}{$\zeta^6$}\  \ \ \Dchaintwo{}{$\zeta^{-1}$}{$\zeta $}{$\zeta^6$} & $\zeta \in G'_{24}$ & $p\not=2,3$
\\
\hline
$13'$ & \Dchaintwo{}{$-\zeta $}{$\zeta^{-1}$}{$-1$}\ \Dchaintwo{}{$1$}{$\zeta$}{$-1$}\ \Dchaintwo{}{$1$}{$\zeta^{-1}$}{\
\ $-\zeta^{2}$}\  \ \ \Dchaintwo{}{$-\zeta^{-1}$}{$-\zeta $}{$-\zeta^2$} & $\zeta \in G'_{8}$ & $p=3$
\\
\hline
14 & \Dchaintwo{}{$\zeta $}{$\zeta^2$}{$-1$}\ \Dchaintwo{}{$-\zeta^{-2}$}{$\zeta^{-2}$}{$-1$} & $\zeta \in G'_5$ &
$p\not=5$
\\
\hline
15 & \Dchaintwo{}{$\zeta $}{$\zeta^{-3}$}{$-1$}\ \Dchaintwo{}{$-\zeta^{-2}$}{$\zeta^3$}{$-1$}\
\Dchaintwo{}{$-\zeta^{-2}$}{$-\zeta^3$}{$-1$}\ \Dchaintwo{}{$-\zeta $}{$-\zeta^{-3}$}{$-1$} & $\zeta \in G'_{20}$ &
$p\not=2,5$
\\
\hline
$ 15'$ & \Dchaintwo{}{$\zeta $}{$\zeta $}{$-1$}\ \Dchaintwo{}{$1$}{$-\zeta $}{$-1$}\ \Dchaintwo{}{$1$}{$\zeta$}{$-1$}\
\Dchaintwo{}{$-\zeta $}{$-\zeta $}{$-1$} & $\zeta \in G'_{4}$ & $p=5$
\\
\hline
16 & \Dchaintwo{}{$\zeta^5 $}{$-\zeta^{-3}$}{$-\zeta$}\ \Dchaintwo{}{$\zeta^5$}{$-\zeta^{-2}$}{$-1$}\
\Dchaintwo{}{$\zeta^3$}{$-\zeta^2$}{$-1$}\ \Dchaintwo{}{$\zeta^3$}{$-\zeta^4$}{$\ -\zeta^{-4}$} & $\zeta \in G'_{15}$ &
$p\not=3,5$
\\
\hline
$16'$ & \Dchaintwo{}{$1$}{$-\zeta^{-1}$}{$-\zeta^2$}\ \Dchaintwo{}{$1$}{$-\zeta$}{$-1$}\
\Dchaintwo{}{$\zeta$}{$-\zeta^{-1}$}{$-1$}\ \Dchaintwo{}{$\zeta$}{$-\zeta^3$}{$\ -\zeta^{-3}$} & $\zeta \in G'_{5}$ &
$p=3$
\\
\hline
$16''$ & \Dchaintwo{}{$\zeta^{-1} $}{$-1$}{$-\zeta$}\ \Dchaintwo{}{$\zeta^{-1}$}{$-\zeta $}{$-1$}\ \Dchaintwo{}{$1$}{$-\zeta^{-1}$}{$-1$}\ \Dchaintwo{}{$1$}{$-\zeta$}{$-\zeta^{-1}$} &
$\zeta \in G'_{3}$ & $p=5$
\\
\hline
17 & \Dchaintwo{}{$-\zeta $}{$-\zeta^{-3}$}{$-1$}\ \Dchaintwo{}{$-\zeta^{-2}$}{$-\zeta^3$}{$-1$} & $\zeta \in G'_7$ &
$p\not=7$
\\
\hline
18 & \Dchaintwo{}{$\zeta^{-1} $}{$-1$}{$-\zeta$} \Dchaintwo{}{$\zeta^{-1}$}{$-\zeta$}{$-1$} \Dchaintwo{}{$1$}{$-\zeta^{-1}$}{$-1$}
 \Dchaintwo{}{$1$}{$-\zeta$}{$-1$} \Dchaintwo{}{$\zeta$}{$-\zeta^{-1}$}{$-1$} \Dchaintwo{}{$\zeta$}{$-1$}{$-\zeta^{-1}\,$} & $\zeta \in G'_3$ & $p=7$
\\
\hline
\end{tabular}
\end{table}}

{\setlength{\unitlength}{1mm} \settowidth{\mpb}{\Dchaintwo{}{$\zeta^{-1} $}{$-1$}{$-\zeta$}
\Dchaintwo{}{$\zeta^{-1}$}{$-\zeta^2$}{$-1$} \Dchaintwo{}{$1$}{$-\zeta^{-1}$}{$-1$},}
\begin{table}[t]\centering
\caption{The exchange graphs of $\mathcal{C}_s(M)$ in Theorem~\ref{Theo:clasi}.}\label{tab.7}
\vspace{1mm}
\begin{tabular}{r|l|l|l|l|l}
\hline
& \text{exchange graphs} &~$n$ &~$l$ &\text{sequences in} $\mathcal{A}^+$ &\text{char} $\Bbbk$
\\
\hline
\hline
1 & \rule[-3\unitlength]{0pt}{8\unitlength}
\begin{picture}
(2,5) \put(1,0){\scriptsize{$\mathcal{D}_{11}$}}
\end{picture}
&1 & 6 & $(0,0)$ &
\\
\hline
2 &
\begin{picture}
(2,5) \put(1,0){\scriptsize{$\mathcal{D}_{21}$}}
\end{picture}
&1 & 4 & $(1,1,1)$ &
\\
\hline
3 &
\begin{picture}
(50,5) \put(1,0){\scriptsize{$\mathcal{D}_{31}$}} \put(6,1){\line(1,0){7}} \put(9,2){\scriptsize{$2$}}
\put(13,0){\scriptsize{$\mathcal{D}_{32}$}} \put(18,1){\line(1,0){7}} \put(22,2){\scriptsize{$1$}} \put(26,0){\scriptsize{$\tau
\mathcal{D}_{31}$}}
\end{picture}
& 3 & 12 & $(1,1,1)$ &
\\
\hline
\rule[-3\unitlength]{0pt}{8\unitlength} 4 &
\begin{picture}
(2,5) \put(1,0){\scriptsize{$\mathcal{D}_{41}$}}
\end{picture}
&1 & 3 & $(2,1,2,1)$ &
\\
\hline
5 &
\begin{picture}
(50,5) \put(1,0){\scriptsize{$\mathcal{D}_{51}$}} \put(6,1){\line(1,0){7}} \put(9,2){\scriptsize{$2$}}
\put(13,0){\scriptsize{$\mathcal{D}_{52}$}}
\end{picture}
& 2 & 6 & $(2,1,2,1)$ &
\\
\hline
6 &
\begin{picture}
(50,5) \put(1,0){\scriptsize{$\mathcal{D}_{61}$}} \put(6,1){\line(1,0){7}} \put(9,2){\scriptsize{$1$}}
\put(13,0){\scriptsize{$\mathcal{D}_{62}$}}
\end{picture}
& 2 & 6 & $(2,1,2,1)$ & $p\not=3$
\\
\hline
\rule[-3\unitlength]{0pt}{8\unitlength} $6' $&
\begin{picture}
(50,5) \put(1,0){\scriptsize{$\mathcal{D}_{6',1}$}} \put(6,1){\line(1,0){7}} \put(9,2){\scriptsize{$1$}}
\put(13,0){\scriptsize{$\mathcal{D}_{6',2}$}}
\end{picture}
& 2 & 6 & $(2,1,2,1)$ & $p=3$
\\
\hline
$6''$ &
\begin{picture}
(50,5) \put(1,0){\scriptsize{$\mathcal{D}_{6'',1}$}}
\end{picture}
& 1 & 3 & $(2,1,2,1)$ &$p\not=2, 3$
\\
\hline
$6'''$ &
\begin{picture}
(50,5) \put(1,0){\scriptsize{$\mathcal{D}_{6''',1}$}}
\end{picture}
& 1 & 3 & $(2,1,2,1)$ &$p=3$
\\
\hline
7 &
\begin{picture}
(50,5) \put(1,0){\scriptsize{$\mathcal{D}_{71}$}} \put(6,1){\line(1,0){7}} \put(9,2){\scriptsize{$2$}}
\put(13,0){\scriptsize{$\mathcal{D}_{72}$}}
\end{picture}
& 2 & 6 & $(2,1,2,1)$ &$p\not=3$
\\
\hline
8 &
 \begin{picture}
(59,15) \put(1,10){\scriptsize{$\mathcal{D}_{81}$}} \put(6,11){\line(1,0){7}} \put(9,12){\scriptsize{$1$}}
\put(13,10){\scriptsize{$\mathcal{D}_{82}$}} \put(2,9){\line(0,-1){5}} \put(0,5){\scriptsize{$2$}} \put(18,11){\line(1,0){7}}
\put(22,12){\scriptsize{$2$}} \put(26,10){\scriptsize{$\mathcal{D}_{83}$}} \put(31,11){\line(1,0){7}}
\put(34,12){\scriptsize{$1$}} \put(38,10){\scriptsize{$\mathcal{D}_{84}$}} \put(43,11){\line(1,0){9}}
\put(47,12){\scriptsize{$2$}} \put(53,10){\scriptsize{$\mathcal{D}_{85}$}} \put(55,9.5){\line(0,-1){5}}
\put(56,5){\scriptsize{$1$}} \put(0,1){\scriptsize{$\tau \mathcal{D}_{85}$}} \put(7,2){\line(1,0){5}}
\put(9,3){\scriptsize{$1$}} \put(13,1){\scriptsize{$\tau \mathcal{D}_{84}$}} \put(20,2){\line(1,0){5}}
\put(22,3){\scriptsize{$2$}} \put(25,1){\scriptsize{$\tau \mathcal{D}_{83}$}} \put(33,2){\line(1,0){6}}
\put(35,3){\scriptsize{$1$}} \put(39,1){\scriptsize{$\tau \mathcal{D}_{82}$}} \put(46,2){\line(1,0){7}}
\put(49,3){\scriptsize{$2$}} \put(53,1){\scriptsize{$\tau \mathcal{D}_{81}$}}
\end{picture}
& 5 & 12 & $(2, 2,1,3,1)$ &$p\not=2,3$
\\
\hline
9 &
\begin{picture}
(60,6) \put(1,0){\scriptsize{$\mathcal{D}_{91}$}} \put(6,1){\line(1,0){7}} \put(9,2){\scriptsize{$2$}}
\put(13,0){\scriptsize{$\mathcal{D}_{92}$}} \put(19,1){\line(1,0){7}} \put(21,2){\scriptsize{$1$}}
\put(26,0){\scriptsize{$\mathcal{D}_{93}$}} \put(32,1){\line(1,0){6}} \put(34,2){\scriptsize{$2$}} \put(38,0){\scriptsize{$\tau
\mathcal{D}_{92}$}} \put(44,1){\line(1,0){6}} \put(47,2){\scriptsize{$1$}} \put(52,0){\scriptsize{$\tau \mathcal{D}_{91}$}}
\end{picture}
& 5 & 12 & $(3,1,2,2,1)$ &$p\not=2,3$
\\
\hline
$9'$ &
\begin{picture}
(60,6) \put(1,0){\scriptsize{$\mathcal{D}_{9',1}$}} \put(7,1){\line(1,0){6}} \put(9,2){\scriptsize{$2$}}
\put(13,0){\scriptsize{$\mathcal{D}_{9',2}$}} \put(19,1){\line(1,0){7}} \put(21,2){\scriptsize{$1$}}
\put(26,0){\scriptsize{$\mathcal{D}_{9',3}$}} \put(32,1){\line(1,0){6}} \put(35,2){\scriptsize{$2$}}
\put(38,0){\scriptsize{$\tau \mathcal{D}_{9',2}$}} \put(44,1){\line(1,0){6}} \put(47,2){\scriptsize{$1$}}
\put(52,0){\scriptsize{$\tau \mathcal{D}_{9',1}$}}
\end{picture}
& 5 & 12 & $(3,1,2,2,1)$ &$p=3$
\\
\hline
10 &
\begin{picture}
(60,6) \put(1,0){\scriptsize{$\mathcal{D}_{10,1}$}} \put(7,1){\line(1,0){7}} \put(10,2){\scriptsize{$2$}}
\put(14,0){\scriptsize{$\mathcal{D}_{10,2}$}} \put(20,1){\line(1,0){7}} \put(22,2){\scriptsize{$1$}}
\put(27,0){\scriptsize{$\mathcal{D}_{10,3}$}}
\end{picture}
& 3 & 6 & $(4,1,2,2,2,1)$ &$p\not=3$
\\
\hline
11 &
\begin{picture}
(60,6) \put(1,0){\scriptsize{$\mathcal{D}_{11,1}$}}
\end{picture}
& 1 & 2 & $(3,1,3,1,3,1)$ &
\\
\hline
12 &
\begin{picture}
(60,6) \put(1,0){\scriptsize{$\mathcal{D}_{12,1}$}} \put(7,1){\line(1,0){7}} \put(10,2){\scriptsize{$2$}}
\put(14,0){\scriptsize{$\mathcal{D}_{12,2}$}} \put(20,1){\line(1,0){7}} \put(22,2){\scriptsize{$1$}}
\put(27,0){\scriptsize{$\mathcal{D}_{12,3}$}}
\end{picture}
& 3 & 6 & $(3,1,3,1,3,1)$ &$p\not=2$
\\
\hline
13 &
\begin{picture}
(60,6) \put(1,0){\scriptsize{$\mathcal{D}_{13,1}$}} \put(7,1){\line(1,0){7}} \put(10,2){\scriptsize{$2$}}
\put(14,0){\scriptsize{$\mathcal{D}_{13,2}$}} \put(20,1){\line(1,0){7}} \put(22,2){\scriptsize{$1$}}
\put(27,0){\scriptsize{$\mathcal{D}_{13,3}$}} \put(33,1){\line(1,0){7}} \put(36,2){\scriptsize{$2$}}
\put(41,0){\scriptsize{$\mathcal{D}_{13,4}$}}
\end{picture}
& 4 & 6 & $(5,1,2,3,1,3,2,1)$ &$p\not=2,3$
\\
\hline
$13'$ &
\begin{picture}
(60,6) \put(1,0){\scriptsize{$\mathcal{D}_{13',1}$}} \put(7,1){\line(1,0){7}} \put(10,2){\scriptsize{$2$}}
\put(14,0){\scriptsize{$\mathcal{D}_{13',2}$}} \put(20,1){\line(1,0){7}} \put(22,2){\scriptsize{$1$}}
\put(27,0){\scriptsize{$\mathcal{D}_{13',3}$}} \put(33,1){\line(1,0){7}} \put(36,2){\scriptsize{$2$}}
\put(41,0){\scriptsize{$\mathcal{D}_{13',4}$}}
\end{picture}
& 4 & 6 & $(5,1,2,3,1,3,2,1)$ &$p=3$
\\
\hline

14 &
\begin{picture}
(60,6) \put(1,0){\scriptsize{$\mathcal{D}_{14,1}$}} \put(7,1){\line(1,0){7}} \put(10,2){\scriptsize{$2$}}
\put(14,0){\scriptsize{$\mathcal{D}_{14,2}$}}
\end{picture}
& 2 & 3 & $(3,1,4,1,3,1,4,1)$ &$p\not=5$
\\
\hline

15 &
\begin{picture}
(60,6) \put(1,0){\scriptsize{$\mathcal{D}_{15,1}$}} \put(7,1){\line(1,0){7}} \put(10,2){\scriptsize{$2$}}
\put(14,0){\scriptsize{$\mathcal{D}_{15,2}$}} \put(20,1){\line(1,0){7}} \put(22,2){\scriptsize{$1$}}
\put(27,0){\scriptsize{$\mathcal{D}_{15,3}$}} \put(33,1){\line(1,0){7}} \put(36,2){\scriptsize{$2$}}
\put(41,0){\scriptsize{$\mathcal{D}_{15,4}$}}
\end{picture}
& 4 & 6 & $(3,1,4,1,3,1,4,1)$ &$p\not=2,5$
\\
\hline
$15' $&
\begin{picture}
(60,6) \put(1,0){\scriptsize{$\mathcal{D}_{15',1}$}} \put(8,1){\line(1,0){7}} \put(11,2){\scriptsize{$2$}}
\put(15,0){\scriptsize{$\mathcal{D}_{15',2}$}} \put(21,1){\line(1,0){7}} \put(23,2){\scriptsize{$1$}}
\put(28,0){\scriptsize{$\mathcal{D}_{15',3}$}} \put(34,1){\line(1,0){7}} \put(37,2){\scriptsize{$2$}}
\put(42,0){\scriptsize{$\mathcal{D}_{15',4}$}}
\end{picture}
& 4 & 6 & $(3,1,4,1,3,1,4,1)$ &$p=5$
\\
\hline
16 &
\begin{picture}
(60,6) \put(1,0){\scriptsize{$\mathcal{D}_{16,1}$}} \put(7,1){\line(1,0){7}} \put(10,2){\scriptsize{$1$}}
\put(14,0){\scriptsize{$\mathcal{D}_{16,2}$}} \put(20,1){\line(1,0){7}} \put(22,2){\scriptsize{$2$}}
\put(27,0){\scriptsize{$\mathcal{D}_{16,3}$}} \put(33,1){\line(1,0){7}} \put(36,2){\scriptsize{$1$}}
\put(41,0){\scriptsize{$\mathcal{D}_{16,4}$}}
\end{picture}
& 4 & 6 & $(2,1,4,1,4,1,2,3)$ & $p\not=3,5$
\\
\hline
$16' $&
\begin{picture}
(60,6) \put(1,0){\scriptsize{$\mathcal{D}_{16',1}$}} \put(8,1){\line(1,0){7}} \put(11,2){\scriptsize{$1$}}
\put(15,0){\scriptsize{$\mathcal{D}_{16',2}$}} \put(22,1){\line(1,0){7}} \put(26,2){\scriptsize{$2$}}
\put(29,0){\scriptsize{$\mathcal{D}_{16',3}$}} \put(36,1){\line(1,0){7}} \put(39,2){\scriptsize{$1$}}
\put(43,0){\scriptsize{$\mathcal{D}_{16',4}$}}
\end{picture}
& 4 & 6 & $(2,1,4,1,4,1,2,3)$ &$p=3$
\\
\hline
$16''$ &
\begin{picture}
(60,6) \put(1,0){\scriptsize{$\mathcal{D}_{16'',1}$}} \put(9,1){\line(1,0){7}} \put(13,2){\scriptsize{$1$}}
\put(17,0){\scriptsize{$\mathcal{D}_{16'',2}$}} \put(25,1){\line(1,0){7}} \put(29,2){\scriptsize{$2$}}
\put(32,0){\scriptsize{$\mathcal{D}_{16'',3}$}} \put(39,1){\line(1,0){7}} \put(42,2){\scriptsize{$1$}}
\put(46,0){\scriptsize{$\mathcal{D}_{16'',4}$}}
\end{picture}
& 4 & 6 & $(2,1,4,1,4,1,2,3)$ &$p=5$
\\
\hline
17 &
\begin{picture}
(50,5) \put(1,0){\scriptsize{$\mathcal{D}_{17,1}$}} \put(9,1){\line(1,0){7}} \put(13,2){\scriptsize{$2$}}
\put(17,0){\scriptsize{$\mathcal{D}_{17,2}$}}
\end{picture}
& 2 & 2 & $(3,1,5,1)^3$ &$p\not=7$
\\
\hline
18 &
\begin{picture}
(78,6) \put(1,0){\scriptsize{$\mathcal{D}_{18,1}$}} \put(8,1){\line(1,0){7}} \put(11,2){\scriptsize{$1$}}
\put(15,0){\scriptsize{$\mathcal{D}_{18,2}$}} \put(22,1){\line(1,0){7}} \put(26,2){\scriptsize{$2$}}
\put(30,0){\scriptsize{$\mathcal{D}_{18,3}$}} \put(37,1){\line(1,0){6}} \put(40,2){\scriptsize{$1$}}
\put(44,0){\scriptsize{$\mathcal{D}_{18,4}$}} \put(50,1){\line(1,0){6}} \put(53,2){\scriptsize{$2$}}
\put(57,0){\scriptsize{$\mathcal{D}_{18,5}$}} \put(63,1){\line(1,0){6}} \put(66,2){\scriptsize{$1$}}
\put(70,0){\scriptsize{$\mathcal{D}_{18,6}$}}
\end{picture}
& 6 & 6 & $(2,1,6,1,2,3)^2$ &$p=7$
\\
\hline
\end{tabular}
\end{table}}

\begin{Remark}
In order to illustrate the exchange graphs of the semi-Cartan graph $\mathcal{C}_s(M)$ in Theorem~\ref{Theo:clasi}, we use the
following notation in Tables~\ref{tab.1}--\ref{tab.7}.
\begin{enumerate}\itemsep=0pt
\item[--] In row~$n$ of Tables~\ref{tab.1}--\ref{tab.6} let $\mathcal{D}_{nl}$ be the~$l$-th Dynkin diagram for all $l\geq1$.
Since the exchange graph of the semi-Cartan graph is labeled, we write $\tau \mathcal{D}_{nl}$ for the graph $\mathcal{D}_{nl}$ where
the two vertices of $\mathcal{D}_{nl}$ change the positions.
\item[--] We also use the notation $(2,1,6,1,2,3)^2=(2,1,6,1,2,3,2,1,6,1,2,3)$ and
$(3,1,5,1)^3=(3,1,5,1,3,1,5,1,3,1,5,1)$ in Table~\ref{tab.7}.
\end{enumerate}
\end{Remark}

\begin{proof}
(1)$\Rightarrow$(3).
Since $\mathcal{B}(V)$ is decomposable and ${\boldsymbol{\Delta}}^{[M]}$ is f\/inite, we obtain that~$M$ admits all ref\/lections and
$\mathcal{R}(M)=(\mathcal{C}(M),(\boldsymbol{\Delta}^{[X]})_{[X]\in \mathcal{X}_2(M)})$ is a~root system of type~$\mathcal{C}(M)$ by~\cite[Corollary~6.12]{HS10}.
Then $\mathcal{W}(M)$ is f\/inite by~\cite[Lemma~5.1]{HS10} since ${\boldsymbol{\Delta}}^{[M]}$ is f\/inite.

(3)$\Rightarrow$(1).
Since~$M$ admits all refections and $\mathcal{W}(M)$ is f\/inite, the set $\boldsymbol{\Delta}^{[M]}$ is f\/inite by~\cite[Lemma 5.1]{HS10}.
Moreover $\mathcal{B}(M)$ is decomposable by~\cite[Corollary~6.16]{a-HeckSchn12a}.

(2)$\Rightarrow$(3).
By Lemma~\ref{le:ifinite} one obtains that~$M$ is~$i$-f\/inite for all $i\in I$.
For $i\in I$, one can determine the Dynkin diagram of $R_i(M)$ by Lemma~\ref{le:Dynkin}.
One observes that it appears in the same row of Table~\ref{tab.6} as $\mathcal{D}$.
Doing the same for all the Dynkin diagrams in the same row of $\mathcal{D}$ implies that~$M$ admits all ref\/lections.
Hence $\mathcal{C}_s(M)$ is well-def\/ined by Proposition~\ref{prop:Cs}.
Now, we identify the objects of $\mathcal{C}_s(M)$ with their Dynkin diagrams.

Assume that $\mathcal{D}$ appears in row~$r$ of one of Tables~\ref{tab.1}--\ref{tab.6}.
Then by the above calculations, the exchange graph of $\mathcal{C}_s(M)$ appears in row~$r$ of Table~\ref{tab.7}.
Then we calculate the smallest integer~$n$ with $(r_2r_1)^n(\mathcal{D})=\mathcal{D}$ and we observe that it appears in the third column
of row~$r$ of Table~\ref{tab.7}.

Then we compute the characteristic sequence $(c_k)_{k\geq1}$ with respect to the f\/irst Dynkin diagram in row~$r$ and the
label~$1$.
We observe that $(c_k)_{k\geq1}$ is the inf\/inite power of the sequence in the f\/ifth column of row~$r$ of
Table~\ref{tab.7}.
Further, we get the numbers $l=6n-\sum\limits_{i=1}^{2n}c_i$.
They appear in the fourth column of Table~\ref{tab.7}.
One checks that $l| 12$ and $(c_1, c_2, \dots,c_{12n/l})\in \mathcal{A}^{+}$ by Corollary~\ref{coro:A+}(2).
The detailed calculations are skipped at this point here.
Then Theorem~\ref{theo.cartan} implies that $\mathcal{C}_s(M)$ is a~f\/inite Cartan graph.
Hence $\mathcal{C}(M)$ is a~f\/inite Cartan graph and $\mathcal{W}(M)$ is a~f\/inite Weyl groupoid by~\cite[Lemma~5.1]{HS10}.

(3)$\Rightarrow$(2).
Since~$M$ admits all ref\/lections, the tuple
\begin{gather*}
\mathcal{C}(M)=\big\{I, \mathcal{X}_2(M), (r_i)_{i\in I}, \big(A^X\big)_{X\in \mathcal{X}_2(M)}\big\}
\end{gather*}
def\/ined in Theorem~\ref{theo.regualrcar} is a~semi-Cartan graph.
In particular,
\begin{gather*}
\mathcal{X}_2(M)=\{[R_{i_1}\cdots R_{i_n}(M)]\in \mathcal{X}_2\,|\, n\in \mathbb{N}_0, \, i_1, \ldots, i_n\in I\}
\end{gather*}
and $R_i(M)=(R_i(M)_j)_{j\in I}$ is def\/ined by equation~\eqref{eq-ri}. 
Moreover, $A^X=(a_{ij}^X)_{i,j\in I}$ for all $X\in \mathcal{X}_2(M)$, where
\begin{gather*}
-a_{ij}^X=\min \big\{m\in \mathbb{N}_0 \,|\, (m+1)_{q_{ii}'}\big({q_{ii}'}^mq_{ij}'q_{ji}'-1\big)=0\big\}
\end{gather*}
by Lemma~\ref{le:aijM} and $(q_{ij}')_{i,j\in I}$ is the braiding matrix of~$X$.
Then $\boldsymbol{\Delta}^{[M]}$ is f\/inite since $\mathcal{W}(M)$ is f\/inite.
Hence all roots are real by~\cite[Proposition~2.12]{c-Heck09b}.
Then $\mathcal{C}(M)$ is a~f\/inite Cartan graph by Theorem~\ref{theo.regualrcar}.
Hence $\mathcal{C}_s(M)$ is a~f\/inite Cartan graph by Proposition~\ref{prop:Cs}.
We can apply Theorem~\ref{theo.cartan} to $\mathcal{C}_s(M)$.

By the implication (2)$\Rightarrow$(3), it is enough to prove that the Dynkin diagram of at least one point in
$\mathcal{C}_s(M)$ is contained in Table~\ref{tab.6}.

Set $X=[M]_s$ and $m=t_{12}^X$.
Let $(c_k)_{k\geq 1}$ be the characteristic sequence of $\mathcal{C}_s(M)$ with respect to~$X$ and the label~$1$.
Then we have $(c_1,c_2,\dots, c_{m})\in \mathcal{A}^+$ and $(c_k)_{k\geq 1}=(c_1,c_2,\dots, c_{m})^{\infty}$~by
Theorem~\ref{theo.cartan}.

If $m=2$ then $(c_1,c_2)=(0,0)$.
Hence $a_{12}^X=a_{21}^X=0$.
Then $q_{12}q_{21}=1$ and $\mathcal{D}=\mathcal{D}_{11}$.

If $m>2$, by Theorem~\ref{Theo:clascar}, one of $(1,1)$, $(1, 2, a)$, $(2, 1, b)$, $(1, 3, 1, b)$ or their transpose,
where $1\leq a\leq 3$ and $3\leq b\leq 5$, is a~subsequence of $(c_1, c_2, \dots, c_{m})$.
Let~$n$ be the smallest integer with $(r_2r_1)^n(X)=X$.
Then $n|m$ by Theorem~\ref{theo.cartan}.
Since $c_{m+k}=c_{k}$ for all $k\in \mathbb{N}$, we have the freedom to assume any position in $(c_k)_{k\geq 1}$, where any of
these subsequences is starting.
Let $c_0=c_m$ and $q_0=q_{12}q_{21}$.

We may assume that $i=1$, $j=2$, but change the labels if necessary.
Next we proceed case by case:

{\samepage Step 1.
If $c_0=c_1=1$, then $a_{12}^X=a_{21}^X=-1$.
Hence $q_0\not=1$.
We distinguish four cases: 1aa, 1ab, 1ba and 1bb.

}

Case 1aa.
If $q_{11}q_{0}=1$ and $q_{22}q_{0}=1$, then $\mathcal{D}=\mathcal{D}_{21}$.

Case 1ab.
If $q_{11}q_{0}=1$, $q_{22}=-1$, and $q_{22}q_{0}\not=1$, then $\mathcal{D}=\mathcal{D}_{31}$.

Case 1ba.
If $q_{11}=-1$, $q_{22}q_{0}=1$, and $q_{11}q_{0}\not=1$, then $\mathcal{D}=\tau \mathcal{D}_{31}$.

Case 1bb.
If $q_{11}=-1$, $q_{22}=-1$, and $q_0\not=-1$, then $\mathcal{D}=\mathcal{D}_{32}$.

\setlength{\unitlength}{1mm} Step 2.
Assume that $(c_0, c_1, c_2)=(1,2,a')$, where $a'\in \{1,2,3\}$.
Then we obtain that $a_{21}^X=-1$, $a_{12}^X=a_{12}^{r_1(X)}=-2$, and $a_{21}^{r_1(X)}=-a'$.
We distinguish four cases: 2aa, 2ab, 2ba and~2bb.

Case 2aa.
If ${q_{11}}^2 q_{0}=1$ and $q_{22}q_{0}=1$, then $\mathcal{D}=\mathcal{D}_{41}$.

Case 2ab.
If ${q_{11}}^2 q_{0}=1$, $q_{22}=-1$, and $q_{22}q_{0}\not=1$, then $\mathcal{D}=\mathcal{D}_{51}$.

Case 2ba.
Assume that $1+q_{11}+{q_{11}}^2=0$, $q_{22}q_{0}=1$, and ${q_{11}}^2 q_{0}\not=1$.
If $p=3$ then $1+q_{11}+{q_{11}}^2=0$ yields ${q_{11}=1}$.
If $q_{22}\not=-1$ then $\mathcal{D}=\mathcal{D}_{6',1}$ and if $q_{22}=-1$ then $\mathcal{D}=\mathcal{D}_{6''',1}$.
Assume that $p\not=3$.
Set $\zeta:=q_{11}$ and $q:=q_{22}$.
Then $q_0=q^{-1}\notin \{1, \zeta^{-1}\}$ since $a_{12}^X=-2$, and $q_0\not=\zeta$ since ${q_{11}}^2 q_{0}\not=1$.
Thus $\mathcal{D}=\mathcal{D}_{61}$ or $\mathcal{D}_{6'',1}$, $p\not=2$.

Case 2bb.
Consider the last case $1+q_{11}+{q_{11}}^2=0$, $q_{22}=-1$ and $q_0\notin \{1,-1,q_{11}, q_{11}^2\}$.

Case 2bba.
If $p=3$ then $q_{11}=1$.
Set $q:=q_{0}$.
By Lemma~\ref{le:Dynkin}, the Dynkin diagrams of $r_1(X)$ and~$X$ are
\begin{gather*}
\Dchaintwo{}{$1$}{$q^{-1}$}{$-q^2$}
\qquad
\Dchaintwo{}{$1$}{$q$}{$-1$}
\end{gather*}
with $q\in \Bbbk^\ast \setminus \{-1,1\}$.
Then $a_{21}^{r_1(X)}\leq -2$ since $(-q^2)q^{-1}=-q\not=1$, and $-q^2\not=-1$.

Case 2bba1.
If $p=3$ and $a'=-a_{21}^{r_1(X)}=2$， then one gets $({-q^2})^2 q^{-1}=1$ or $1+(-q^2)+(- q^2)^2=0$.
If $({-q^2})^2 q^{-1}=1$ then $q=1$, which is a~contradiction.
Hence $-q^2=1$ from the second equation since $p=3$.
Then $\mathcal{D}=\mathcal{D}_{9',2}$.

Case 2bba2.
If $p=3$ and $a'=-a_{21}^{r_1(X)}=3$, then one has $({-q^2})^3 q^{-1}=1$ or $1+(- q^2)+(- q^2)^2+(-q^2)^3=0$.
The f\/irst equation $({-q^2})^3 q^{-1}=1$ yields $(-q)^5=1$, hence $\mathcal{D}=\mathcal{D}_{16',2}$.
If $1+(- q^2)+(- q^2)^2+(-q^2)^3=0$, then $(1-q^2)(1+q^4)=0$ and hence $q\in G'_8$.
Then $\mathcal{D}=\mathcal{D}_{13',2}$.

Case 2bbb.
We now suppose that $p\not=3$.
Set $\zeta:=q_{11}$ and $q:=q_{0}$.
Hence the Dynkin diagram of $r_1(X)$ is \Dchaintwo{}{$\zeta$}{${(\zeta q)}^{-1}$}{$
\ \ \ -\zeta q^2$} \ \ \ \ with $\zeta \in G'_3$, $q\in \Bbbk^*\setminus \{1, -1, \zeta, \zeta^{-1}\}$.
Since $a'\in \{1, 2, 3\}$, we distinguish three cases:
\begin{enumerate}\itemsep=0pt
\item[(b1)] $p\not=3$, $a_{21}^{r_1(X)}=1$,
\item[(b2)] $p\not=3$, $a_{21}^{r_1(X)}=2$,
\item[(b3)] $p\not=3$, $a_{21}^{r_1(X)}=3$.
\end{enumerate}

Case 2bbb1.
If the condition (b1) holds, then one gets $(-\zeta q^2)({\zeta q})^{-1}=1$ or $1+(-\zeta q^2)=0$.
If $(-\zeta q^2)({\zeta q})^{-1}=1$, then $q=-1$, which is a~contradiction.
If $1+(-\zeta q^2)=0$ then $\zeta^2=q^2$ and hence $q=-\zeta$.
Then $\mathcal{D}=\mathcal{D}_{71}$.

Case 2bbb2.
If the condition $(b2)$ holds, then $({-\zeta q^2})^2 (\zeta q)^{-1}=1$ or $\sum\limits_{i=0}^2(-\zeta q^2)^i=0$.

Case 2bbb2a.
Consider the equation $({-\zeta q^2})^2 (\zeta q)^{-1}=1$.
Then $\zeta q^3=1$ and hence $q\in G'_9$ since $\zeta\in G'_3$ and $p\not=3$.
Hence $\mathcal{D}=\mathcal{D}_{10,2}$.

Case 2bbb2b.
If $\sum\limits_{i=0}^2(-\zeta q^2)^i=0$, then $-\zeta q^2\in \{\zeta, \zeta^{-1}\}$.
Hence $q^2=-1$ or $-q^2=\zeta$.

Case 2bbb2b1.
If $q^2=-1$, then $p\not=2$ and the Dynkin diagram of~$X$ is$\!\!\!$ \Dchaintwo{}{$\zeta $}{$q$}{$-$1} \ with $q\in G'_4$, $\zeta
\in G'_3$.
Set $\eta:={\zeta}^2 q^{-1}$.
Then $\eta \in G'_{12}$, $\zeta=-\eta^2$, and $q=\eta^3$.
Hence $\mathcal{D}=\mathcal{D}_{92}$.

Case 2bbb2b2.
If $-q^2=\zeta$, then $q\in G'_{12}$ and $p\not=2$ since $q\not=\zeta^{-1}$.
Hence $\mathcal{D}=\mathcal{D}_{81}$.

Case 2bbb3.
If the condition $(b3)$ holds, then $({-\zeta q^2})^3 (\zeta q)^{-1}\!=1$ or $\sum\limits_{i=0}^3(-\zeta q^2)^i\!=0$,
$1-\zeta q^2\not=0$.

Case 2bbb3a.
Consider the equation $({-\zeta q^2})^3 (\zeta q)^{-1}=1$, that is $-q^5=\zeta$.
Hence $q=-\zeta^{-1}$ or $-q\in G_{15}'$, $p\not=5$.

Case 2bbb3a1.
If $-q\in G'_{15}$, $p\not=3,5$, and $\zeta=-q^5$, then $\mathcal{D}=\mathcal{D}_{16,2}$.

Case 2bbb3a2.
If $q=-\zeta^{-1}$ then $p\not=2$ since $q\not=\zeta^{-1}$.
Hence the Dynkin diagrams of $r_1(X)$,~$X$ and $r_2(X)$, respectively, are
\begin{gather*}
\Dchaintwo{}{$\zeta$}{$-1$}{$-\zeta^{-1}$}
\qquad
\Dchaintwo{}{$\zeta$}{$-\zeta^{-1}$}{$-1$}
\qquad
\Dchaintwo{}{$1$}{$-\zeta$}{$-1$}
\qquad
\end{gather*}
with $\zeta\in G'_3$.
Then one obtains that $a_{12}^{r_2(X)}=1-p$.
We distinguish four cases.

Case 2bbb3a2a.
If $p=5$, then $\mathcal{D}=\mathcal{D}_{16''2}$.

Case 2bbb3a2b.
If $p=7$, then $\mathcal{D}=\mathcal{D}_{18,2}$.

Case 2bbb3a2c.
If $p=6s+1$ $(s \geq 2)$, then the Dynkin diagrams of $r_1(X)$,~$X$, $r_2(X)$, $r_1r_2(X)$, $r_2r_1r_2(X)$, and
$(r_1r_2)^2(X)$, respectively, are
\begin{gather*}
\Dchaintwo{}{$\zeta$}{$-1$}{$-\zeta^{-1}$}
\qquad
\Dchaintwo{}{$\zeta$}{$-\zeta^{-1}$}{$-1$}
\qquad
\Dchaintwo{}{$1$}{$-\zeta$}{$-1$}
\qquad
\Dchaintwo{}{$1$}{$-\zeta^{-1}$}{$-1$}
\qquad
\Dchaintwo{}{$\zeta^{-1}$}{$-\zeta$}{$-1$}
\qquad
\Dchaintwo{}{$\zeta^{-1}$}{$-1$}{$-\zeta$}
\end{gather*}
with $\zeta\in G'_3$.
Hence $n=6$ in Theorem~\ref{theo.cartan} and $(c_k)_{k\geq 0}=(2, 3, 2, 1, p-1, 1, 2, 3, 2, 1, p-1, 1)^{\infty}$.
Then $l=20-2p<0$, which is a~contradiction to Theorem~\ref{theo.cartan}.

Case 2bbb3a2d.
If $p=6s+5$, where $s\geq 1$, then the Dynkin diagrams of $r_1(X)$,~$X$, $r_2(X)$ and $r_1r_2(X)$, respectively, are
\setlength{\unitlength}{1mm}
\begin{gather*}
\Dchaintwo{}{$\zeta$}{$-1$}{$-\zeta^{-1}$}
\qquad
\Dchaintwo{}{$\zeta$}{$-\zeta^{-1}$}{$-1$}
\qquad
\Dchaintwo{}{$1$}{$-\zeta$}{$-1$}
\qquad
\Dchaintwo{}{$1$}{$-\zeta^{-1}$}{$-\zeta$}
\qquad
\end{gather*}
with $\zeta\in G'_3$.
Then $n=4$ and $(c_{k})_{k \geq 0}=(2, 3, 2, 1, p-1, 1, p-1, 1)^{\infty}$.
Hence $l=16-2p<0$.
Again, one gets a~contradiction.

Case 2bbb3b.
Consider the equation $0=\sum\limits_{i=0}^3(-\zeta q^2)^i=(1-\zeta q^2)(1+\zeta^2 q^4)$, where $\zeta q^2\not=1$.
One gets $\zeta=-q^{4}$.
If $p=2$, then $\zeta^4=q^{4}$ and hence $\zeta=q$, which is a~contradiction to $\zeta q^2\not=1$.
Otherwise $q\in G'_{24}$ and $\mathcal{D}=\mathcal{D}_{13,2}$.

Step~3.
Now we change the label.
It means that $(c_k)_{k\geq 1}$ is the characteristic sequence of~$\mathcal{C}_s(M)$ with respect to~$X$ and the label~$2$.

Assume that $(c_0, c_1, c_2)=(2,1,b')$, where $b'\in \{3,4,5\}$.
Then we obtain that $a_{12}^X=-2$, $a_{21}^X=-1$ and $a_{12}^{r_2(X)}=-b'$.
If ${q_{11}}^2 q_{0}=1$ and $q_{22}=-1$ then $a_{12}^{r_2(X)}=-2$, which is a~contradiction.
If $q_0q_{22}=1$ then $a_{12}^{r_2(X)}=a_{12}^X=-2$, which is again a~contradiction.
Suppose now that $1+q_{11}+{q_{11}}^2=0$, $q_{22}=-1$ and $q_0\notin \{1, -1, q_{11}^{-2}\}$. Since $a_{12}^X=-2$, we
also obtain that $q_0\not=q_{11}^{-1}$.

Case 3a.
If $p=3$, then by setting $q:=q_{0}$, the Dynkin diagrams of~$X$ and $r_2(X)$, respectively, are
\begin{gather}
\label{diam1}
\Dchaintwo{}{$1$}{$q$}{$-1$}
\qquad
\Dchaintwo{}{$-q$}{$q^{-1}$}{$-1$}
\qquad
q\in \Bbbk^\ast \setminus \{-1,1\}.
\end{gather}

Case 3a1.
If $p=3$ and $b'=3$, then one gets $({-q})^3 q^{-1}=1$ or $\sum\limits_{i=0}^3(- q)^i=0$.
Since $q\not=1$, both equations imply that $q^2=-1$.
Hence $\mathcal{D}=\mathcal{D}_{9',2}$.

Case 3a2.
If $p=3$ and $b'=4$, then one has $({-q})^4 q^{-1}=1$ or $\sum\limits_{i=0}^4(- q)^i=0$.
If $({-q})^4 q^{-1}=1$ then $q=1$, which is a~contradiction to~\eqref{diam1}.
If $\sum\limits_{i=0}^4(- q)^i=0$ then $q^5=-1$ and $\mathcal{D}=\mathcal{D}_{16',2}$.

Case 3a3.
If $p=3$ and $b'=5$, then one obtains $({-q})^5 q^{-1}=1$ or $\sum\limits_{i=0}^5(- q)^i=0$.
If $({-q})^5 q^{-1}=1$ then $q\in G'_8$ and $\mathcal{D}=\mathcal{D}_{13',2}$.
Consider the equation $0=\sum\limits_{i=0}^5(- q)^i=(1-q)(1+q^2+q^4)$.
Then $q^2=1$ since $p=3$ and $q\not=1$, which is a~contradiction to~\eqref{diam1}.

Case 3b.
We now consider the cases in which the condition $p\not=3$ holds.
Set $\zeta:=q_{11}$ and $q:=q_{0}$.
The Dynkin diagram of $r_2(X)$ is
\begin{gather}
\label{diam2}
\Dchaintwo{}{$-\zeta q$}{$q^{-1}$}{$-1$}
\qquad
\zeta \in G'_3,
\qquad
q\in \Bbbk^*\setminus \big\{1, -1, \zeta, \zeta^{-1}\big\}.
\end{gather}

Case 3b1.
If $p\not=3$ and $b'=3$, then one gets ${(-\zeta q)}^3{q}^{-1}=1$ or $\sum\limits_{i=0}^3(-\zeta q)^i=0$.
If ${(-\zeta q)}^3{q}^{-1}=1$ then $q\in G'_4$, $p\not=2$ and $\mathcal{D}=\mathcal{D}_{92}$.
If $0=\sum\limits_{i=0}^3(-\zeta q)^i=(1-\zeta q)(1+(\zeta q)^2)$, then $(\zeta q)^2=-1$ and $p\not=2$ since
$q\not=\zeta^{-1}$.
Hence $\mathcal{D}=\mathcal{D}_{81}$.

Case 3b2.
If $p\not=3$ and $b'=4$, then one gets ${(-\zeta q)}^4{q}^{-1}=1$ or $\sum\limits_{i=0}^4(-\zeta q)^i=0$.
If ${(-\zeta q)}^4{q}^{-1}=1$ then $\zeta=q^{-3}$.
Since $q\notin G'_3$, one obtains $q\in G'_9$ and $\mathcal{D}=\mathcal{D}_{10,2}$.
The equation $\sum\limits_{i=0}^4(-\zeta q)^i=0$ gives $\zeta=-q^{5}$.
Since $\zeta\in G'_3$, one gets $-q\in G'_3$, $\zeta=-q^{-1}$, $p=5$ or $-q\in G'_{15}$, $p\not=3,5$.
If $-q\in G'_3$ then $\mathcal{D}=\mathcal{D}_{16'',2}$ and if $-q\in G'_{15}$ then $\mathcal{D}=\mathcal{D}_{16,2}$.

Case 3b3.
If $p\not=3$ and $b'=5$, then one gets ${(-\zeta q)}^5{q}^{-1}=1$ or $\sum\limits_{i=0}^5(-\zeta q)^i=0$, $-\zeta
q\not=1$.
If ${(-\zeta q)}^5{q}^{-1}=1$ and $p=2$ then $q\zeta^{-1}=1$, which is a~contradiction to~\eqref{diam2}.
If ${(-\zeta q)}^5{q}^{-1}=1$ and $p\not=2$ then $\zeta=-q^4$ and $\mathcal{D}=\mathcal{D}_{13,2}$.
Consider $0=\sum\limits_{i=0}^5(-\zeta q)^i=(1-\zeta q)(1+(-\zeta q)^2+(-\zeta q)^4)$.
Since $q\notin \{1,-1,\zeta,\zeta^{-1}\}$, one gets $p\not=2$ and $q^3=-1$.
Since $-\zeta q\not=1$, one gets $q=-\zeta$ and $a_{12}^{r_2(X)}=-2$, which is a~contradiction.

Step 4.
Again we use the same labeling as in steps $1$ and~$2$.
Assume that $(c_0, c_1, c_2, c_3)=(1, 3, 1, c')$, where $c'\in \{3, 4, 5\}$.
Then $a_{21}^X=-1$ and $a_{12}^X=-3$.
We distinguish four cases: 4aa, 4ab, 4ba and 4bb.

Case 4aa.
If ${q_{11}}^3 q_{0}=1$ and $q_{22}q_{0}=1$, then $\mathcal{D}=\mathcal{D}_{11,1}$.

Case 4ab.
Set $q:=q_{11}$.
If ${q_{11}}^3 q_{0}=1$ and $q_{22}=-1$, then the Dynkin diagrams of~$X$, $r_1(X)$ and $r_2r_1(X)$, respectively, are
\begin{gather*}
\Dchaintwo{}{$q$}{$q^{-3}$}{$-1$}
\qquad
\Dchaintwo{}{$q$}{$q^{-3}$}{$-1$}
\qquad
\Dchaintwo{}{-${q}^{-2}$}{$q^3$}{$-1$}
\qquad
q\in \Bbbk^*\setminus \{1,-1\},
\qquad
q\notin G'_3.
\end{gather*}
Then $-a_{12}^{r_2r_1(X)}=c'\in \{3,4,5\}$.
Hence we distinguish three cases: 4ab1, 4ab2 and 4ab3.

Case 4ab1.
If $c'=3$, then $(-{q}^{-2})^3 {q}^3=1$ or $\sum\limits_{i=0}^3(-{q}^{-2})^i=0$, $1-q^{-2}\not=0$.
If $(-{q}^{-2})^3 {q}^3=1$ then $\mathcal{D}=\mathcal{D}_{11,1}$, where $q^3=-1$.
If $\sum\limits_{i=0}^3(-{q}^{-2})^i=0$, then $\mathcal{D}=\mathcal{D}_{12,1}$.

Case 4ab2.
If $c'=4$, then $(-{q}^{-2})^4 {q}^3=1$ or $\sum\limits_{i=0}^4(-{q}^{-2})^i=0$.

Case 4ab2a.
The equation $(-{q}^{-2})^4 {q}^3=1$ gives $q^5=1$ and $p\not=5$, since $q\not=1$.
Hence $\mathcal{D}=\mathcal{D}_{14,1}$.

Case 4ab2b.
Consider the equation $\sum\limits_{i=0}^4(-{q}^{-2})^i=0$.
One gets $q^{10}=-1$.
If $p=2$ then $\mathcal{D}=\mathcal{D}_{14,1}$.
If $p=5$ then $q^2=-1$ and $\mathcal{D}=\mathcal{D}_{15',1}$.
If $p\not=2, 5$, then $q\in G'_{20}$ and $\mathcal{D}=\mathcal{D}_{15,1}$.

Case 4ab3.
If $c'=5$, then $(-{q}^{-2})^5 {q}^3=1$ or $\sum\limits_{i=0}^5(-{q}^{-2})^i=0$.

Case 4ab3a.
Consider the equation $(-{q}^{-2})^5 {q}^3=1$, which gives ${-q}^7=1$.
Since $q\not=-1$, one gets $p\not=7$ and $-q\in G'_{7}$.
Hence $\mathcal{D}=\mathcal{D}_{17,1}$.

Case 4ab3b.
Consider the equation $0=\sum\limits_{i=0}^5(-{q}^{-2})^i=(1-q^{-2})(1+q^{-4}+q^{-8})$.
Since $q^2\not=1$, one gets $1+q^{-4}+q^{-8}=0$.
If $p=3$ then $q\in G'_4$ since $q^2\not=1$.
Hence $a_{12}^{r_2r_1(X)}=-2$, which is a~contradiction.
Then $p\not=3$ and $q^{-4}\in G'_3$.
Since $c'=5$, one gets $q\in G'_6$ or $q\in G'_{12}$.
If $q\in G'_6$, then $a_{12}^{r_2r_1(X)}=-3$, which is a~contradiction.
If $q\in G'_{12}$, then $a_{12}^{r_2r_1(X)}=-2$, which is again a~contradiction.

Case 4ba.
The conditions $1+q_{11}+{q_{11}}^2+{q_{11}}^3=0$, $q_{11}\not=-1$ and $q_{22}q_{0}=1$ hold.
Then $q_{11}\in G'_4$ and $p\not=2$ since $a_{12}^X=-3$. Set $\zeta:=q_{11}$ and $q=q_{22}$.
The Dynkin diagram of~$r_1(X)$ is \Dchaintwo{}{$\zeta$}{$-q$}{$
\qquad
\zeta q^{-2}$} \ \ \ with $\zeta \in G'_4$ and $q\in \Bbbk^*\setminus \{1, -1, \zeta, \zeta^{-1}\}$.
Since $-a_{21}^{r_1(X)}=c_2=1$, one gets $(\zeta q^{-2})(-q)=1$ or $\zeta q^{-2}=-1$.
If $(\zeta q^{-2})(-q)=1$ then $q=-\zeta$, which is a~contradiction.
If $\zeta q^{-2}=-1$ then $q\in G'_8$ and $\mathcal{D}=\mathcal{D}_{12,3}$.

Case 4bb.
Consider the last case: $1+q_{11}+q_{11}^2+q_{11}^3=0$, $q_{22}=-1$ and $q_{11}\not=-1$.
Then $q_{11}^2=-1$ and $p\not=2$.
Set $\zeta:=q_{11}$ and $q=q_{0}$.
The Dynkin diagram of $r_1(X)$ is$\!\!\!$ \Dchaintwo{}{$\zeta$}{$\!-q^{-1}$}{$\,\,\,-\zeta q^3$} \ \ with $q\in \Bbbk^*\setminus \{1, -1,
\zeta, \zeta^{-1}\}$ and $\zeta \in G'_4$. Since $-a_{21}^{r_1(X)}=c_2=1$, one has $(-q^{-1})(-\zeta q^3)=1$ or $\zeta q^3=1$.

Case 4bb1.
If $\zeta q^3=1$, then $\zeta=q^{-3} \in G'_4$.
If $p=3$, then $\zeta=q$, which is a~contradiction.
Hence $p\not=3$ and $\mathcal{D}=\tau\mathcal{D}_{85}$.

Case 4bb2.
The condition $(-q^{-1})(-\zeta q^3)=1$ holds.
Then $\zeta=q^{-2}\in G'_4$.
Hence $q\in G'_8$ and $\mathcal{D}=\mathcal{D}_{12,2}$.

By checking all cases in Theorem~\ref{Theo:clascar}, the proof of Theorem~\ref{Theo:clasi} is completed.
\end{proof}

\begin{Remark}
Assume that char~$\Bbbk=p>0$.
Let~$V$ be a~two-dimensional braided vector space of diagonal type.
Let $(x_1, x_2)$ be a~basis of~$V$ and $(q_{ij})_{1\leq i, j\leq2}\in (\Bbbk\setminus\{0\})^{2\times 2}$ satisfying
\begin{gather*}
c(x_i\otimes x_j)=q_{ij}x_j\otimes x_i
\end{gather*}
for any $i$, $j$.
\begin{itemize}\itemsep=0pt
\item[--] By~\cite[Corollary~6]{a-Heck04e}, $\dim_\Bbbk\mathcal{B}(V)<\infty$ if and only if $M=(\Bbbk x_1,\Bbbk x_2)$ admits
all ref\/lections, $\mathcal{W}(M)$ is f\/inite, and for all points $\mathcal{D}$ of $\mathcal{C}_s(M)$, the labels of the vertices of $\mathcal{D}$
are roots of unity (including~$1$).
Therefore with Theorem~\ref{Theo:clasi} one can easily decide whether $\dim_\Bbbk \mathcal{B}(V)$ is f\/inite.
\item[--] If the Dynkin diagram of~$V$ appears in the row $18$ of Table~\ref{tab.6}, then the Weyl groupoid of $\mathcal{B}(V)$
is not appearing for Nichols algebras in characteristic zero.
\end{itemize}
\end{Remark}

\subsection*{Acknowledgements}

It is a~pleasure to thank N.~Andruskiewitsch and H.-J.~Schneider for a~very fruitful discussion on some details of this
topic.
The authors thank the referees for their helpful comments and suggestions.
J.~Wang is supported by China Scholarship Council.

\pdfbookmark[1]{References}{ref}
\LastPageEnding

\end{document}